\documentclass[10pt,a4paper]{article}
\usepackage[a4paper]{geometry}
\usepackage{amssymb,latexsym,amsmath,amsfonts,amsthm}
\usepackage{graphicx,subfigure}
\usepackage{epsfig}

\advance\textheight by \topskip

\newtheorem{teo}{\sc Theorem}[section]

\newtheorem{cor}{\sc Corollary}[section]
\newtheorem{lemma}{\sc Lemma}[section]
\newtheorem{pro}{\sc Proposition}[section]

\theoremstyle{definition}
\newtheorem{defi}{\sc Definition}[section]

\theoremstyle{remark}
\newtheorem{rem}[teo]{\sc Remark}

\newtheorem{notation}{\sc Notation}

\newcommand{\bte}{\begin{teo}}
\newcommand{\ete}{\end{teo}}
\newcommand{\bc}{\begin{cor}}
\newcommand{\ec}{\end{cor}}
\newcommand{\bp}{\begin{pro}}
\newcommand{\ep}{\end{pro}}
\newcommand{\bl}{\begin{lemma}}
\newcommand{\el}{\end{lemma}}
\newcommand{\bd}{\begin{defi}}
\newcommand{\ed}{\end{defi}}
\newcommand{\bno}{\begin{notation}}
\newcommand{\eno}{\end{notation}}

\newcommand{\bca}{\begin{cases}}
\newcommand{\eca}{\end{cases}}

\newcommand{\bq}{\begin{equation}}
\newcommand{\eq}{\end{equation}}
\newcommand{\btabu}{\begin{table}}
\newcommand{\etabu}{\end{table}}
\newcommand{\bt}{\begin{tabular}}
\newcommand{\et}{\end{tabular}}
\newcommand{\ba}{\begin{array}}
\newcommand{\ea}{\end{array}}
\newcommand{\br}{\begin{eqnarray}}
\newcommand{\er}{\end{eqnarray}}
\newcommand{\brn}{\begin{eqnarray*}}
\newcommand{\ern}{\end{eqnarray*}}
\newcommand{\benu}{\begin{enumerate}}
\newcommand{\eenu}{\end{enumerate}}
\newcommand{\bite}{\begin{itemize}}
\newcommand{\eite}{\end{itemize}}

\newcommand{\supp}{\operatorname{supp }}

\title{High order recurrence relation, Hermite-Pad\'e approximation, and Nikishin systems}

\author{D. Barrios Rolan\'{\i}a\footnotemark[1] , J. S. Geronimo\footnotemark[2], G. L\'opez Lagomasino\footnotemark[3]}

\begin{document}

\maketitle

\renewcommand{\thefootnote}{\fnsymbol{footnote}}

\footnotetext[1]{Ingenier\'{\i}a Civil: Hidr\'aulica y Ordenaci\'on del Territorio, Universidad Polit\'ecnica de Madrid, Madrid, Spain. email: dolores.barrios.rolania\symbol{'100}upm.es}
\footnotetext[2]{Department of Mathematics, Georgia Institute of Technology, Atlanta, USA.   email:  geronimo\symbol{'100}math.gatech.edu}
\footnotetext[3]{Departamento de Matem\'aticas, Universidad Carlos III de Madrid, Avda. Universidad 30, 28911 Legan\'es, Madrid, Spain. email:
lago\symbol{'100}math.uc3m.es.\\  D.B.R. was partially supported by research grant MTM2014-54053-P of Ministerio de Econom\'{\i}a y Competitividad, Spain, J.S.G. was partially supported by Simon's Foundation Collaboration Grant and expresses his gratitude to the Department of Mathematics of the University Carlos III de Madrid for its hospitality during the semester spent there making use of the grant C\'atedra de Excelencia Universidad Carlos III de Madrid-Banco Santander, and G.L.L. was supported by research grant MTM2015-65888-C4-2 of Ministerio de Econom\'{\i}a y Competitividad, Spain.}

\begin{abstract} The study of sequences of polynomials satisfying high order recurrence relations is connected with the asymptotic behavior of multiple orthogonal polynomials, the convergence properties of type II Hermite-Pad\'e approximation, and eigenvalue distribution of banded Toeplitz matrices. We present some  results for the case of recurrences  with constant coefficients which match what is known for the Chebyshev polynomials of the first kind. In particular, under appropriate assumptions, we show that  the sequence of polynomials satisfies multiple orthogonality relations with respect to a Nikishin system of measures. \end{abstract} \vspace{1cm}


{\it A.M.S. Subject Classification.} Primary: 30E10, 42C05;
Secondary: 41A20.

\section{Introduction}

Let us consider the sequence $(Q_n)_{n \geq 0}$ of monic polynomials defined using a general $p+2$ term recurrence relation with constant coefficients
\begin{equation}
\label{rec}
\lambda Q_{n}(\lambda) =  Q_{n+1}(\lambda) + a_0 Q_n (\lambda) +  \cdots + a_p Q_{n-p}(\lambda), \quad a_j \in {\mathbb{C}}  \quad j=0,\ldots,p,\quad a_p  \neq 0,
\end{equation}
with initial conditions $Q_0(\lambda) \equiv 1, Q_{-1}(\lambda) = \cdots = Q_{-p}(\lambda) \equiv 0$.

Define the infinite banded Hessenberg matrix
\begin{equation} \label{eq:A}
A := \left(
\begin{array}{cccc}
a_0  & 1 & 0 & \cdots \\
a_1 & a_0 & 1 & \ddots\\
\vdots & \ddots & \ddots & \ddots \\
a_p & a_{p-1} & \ddots & \ddots \\
0 & a_p & a_{p-1} & \ddots \\
\vdots & \ddots & \ddots & \ddots
\end{array}
\right)
\end{equation}
If ${\mathbf{Q}} = (Q_0,\ldots, Q_n, \ldots)^t$, where $(\cdot)^t$ means taking transpose, the recurrence relation can be viewed  in matrix form as
\[ A \mathbf{Q} = \lambda \mathbf{Q}.
\]
It is well known and easy to verify that the zeros of $Q_n$ are the eigenvalues associated with the $n$th principal section $A_n$ of $A$.

For $n\in{\mathbb{Z}}_{+}$ we define the multi-index
\begin{equation}\label{multi:index}
{\mathbf{n}}:=(n_1,\ldots,n_p):=(\underset{\textrm{$k$ times}}{\underbrace{m+1,\ldots,m+1}},\underset{\textrm{$p-k$
times}}{\underbrace{m,\ldots,m}})\in{\mathbb{ Z}}^p_{+},
\end{equation}
where $m\in{\mathbb{Z}}_+$ and $k\in\{0,1,\ldots,p-1\}$ are such that $n=mp+k$. Note that we have
$|{\mathbf{n}}|:=n_{1}+n_{2}+\cdots+n_{p}=n$.

The study of the asymptotic behavior of sequences of polynomials $(Q_n)_{ n \geq 0}$ satisfying higher order recurrence relations with constant or periodic coefficients and their connections with Hermite-Pad\'e approximation has been a subject of major interest. To our knowledge, the first paper in this direction for four term recurrences is \cite{kn:Apt}. The case when $a_0 =\cdots = a_{p-1} =0$,  related with Faber polynomials associated with hypocycloidal domains,  was considered in \cite{HS} and \cite{AKS}. The situation when $a_0 =\cdots = a_{p-1} =0$ and the coefficient corresponding to $a_p$ is periodic was analyzed in \cite{LopezGarcia} and \cite{DL}. A situation in which all the coefficients have period $p$, which extends the notion of Chebyshev polynomials to the context of multiple orthogonality, was treated in \cite{DLL} (see also \cite{LopRoc}). These questions are linked as well with the study of vector valued continued fractions \cite{kn:Par} and \cite{VIse}, and the spectral analysis of banded Hessenberg and Toeplitz operators \cite{kn:Kal} and \cite{kn:DK}.

For every $\lambda \in {\mathbb{C}}$ let $z_j(\lambda), j=0,\ldots,p,$ denote the $p+1$ solutions of the algebraic equation
\begin{equation}
\label{alg2}
  a(z) - \lambda = \sum_{k=-1}^p a_kz^k - \lambda = 0, \qquad a_{-1} = 1.
\end{equation}
The rational function $a(z)$ is called the symbol of the Toeplitz matrix $A$. Its inverse $z(\lambda)$ establishes a one to one correspondence between a $p+1$ sheeted Riemann surface $\mathcal{R}$ of genus zero (see  Hurwitz' formula \cite[p. 52]{kn:Mir}) and the extended complex plane $\overline{\mathbb{C}}$.
We order the solutions of the algebraic equation by absolute value so that
\begin{equation}
\label{order}
  0 < |z_0(\lambda)| \leq \cdots \leq |z_{p}(\lambda)|, \qquad \lambda \in \mathbb{C}.
\end{equation}
When all inequalities are strict the order is unambiguously defined. If equalities occur, choose an arbitrary ordering so that the inequalities remain true.
Set
\[ \Gamma_k = \{\lambda \in {\mathbb{C}}\,:\,|z_{k-1}(\lambda)| = |z_{k}(\lambda)|\}, \qquad k=1,\ldots,p.
\]
It is known that each $\Gamma_k$ is a disjoint union of finitely many open analytic arcs and a finite number of exceptional points with no isolated points. All $\Gamma_k$ are unbounded except $\Gamma_1$ which is compact and connected (see \cite[Proposition 3.2]{kn:DK} and \cite[Theorems 11.9, 11.19]{BG}). We call attention to the fact that the labeling of the $z_k(\lambda)$ and $\Gamma_k$ in this paper does not match the one adopted in \cite{kn:DK}.
The point $\lambda = \infty$ is a branch point of order $p-1$ of the algebraic function and it is also a simple zero (so $z_0(\infty) = 0$).
Set
\[ {\mathcal{R}}_k = \{ \lambda \in {\mathbb{C}}: |z_{k-1}(\lambda)| < |z_{k}(\lambda)| < |z_{k+1}(\lambda)|\}, \qquad k=0,\ldots,p
\]
with the convention $z_{-1}(\lambda) \equiv 0, z_{p+1}(\lambda) \equiv \infty$.
Then
\[{\mathcal{R}} = \overline{\cup_{k=0}^p {\mathcal{R}}_k}.\]

Let ${\mathbf e}_j$ be the unit vector in $\ell^2$ with $1$ in the $j+1$ component. Define
\begin{equation}\label{resolvf}
g_{j}(\lambda):=(R_{\lambda} {\mathbf e}_{j-1},{\mathbf e}_{0}),\qquad 1\leq j\leq p,
\end{equation}
where $R_{\lambda}=(\lambda I-A)^{-1}$ is the resolvent operator, and $(\cdot,\cdot)$ is the standard inner product in $\ell^{2}$.
For each $j=1,\ldots,p$, we introduce a linear functional $L_{j}$ defined in the space of polynomials by
\begin{equation} \label{L}
L_{j}(z^{n})=(A^{n} {\mathbf v}_{j},{\mathbf e}_{0}),
\end{equation}
where ${\mathbf v}_{j}:={\mathbf e}_{0}+\cdots+{\mathbf e}_{j-1}$. Notice that $A$ defines a bounded operator and thus the functions $g_j$ admit convergent Laurent expansions in a neighborhood of $\infty$, and $g_j(\infty) = 0, j=1,\ldots,p$.

\begin{teo} \label{teo:1} We have:
\begin{itemize}
\item[i)] $(g_1(\lambda),\ldots,g_p(\lambda)) = (z_0(\lambda), \ldots,z_0^p(\lambda))$. More precisely, both vector functions have the same Laurent expansions at $\infty$ componentwise.
\item[ii)] For each $j\in\{1,\ldots,p\}$, the polynomial $Q_{n}$ satisfies the following orthogonality conditions
\[
L_{j}(\lambda^{k} Q_{n}(\lambda) )=0,\qquad k=0,\ldots,n_{j}-1.
\]
\item[iii)]
The vector rational function
\[
\Big(\frac{Q_{n-1}}{Q_{n}},\frac{Q_{n-2} }{Q_{n}},\ldots,\frac{Q_{n-p} }{Q_{n}}\Big)
\]
is a type II Hermite-Pad\'e approximant to the system $(z_0,z_0^2,\ldots,z_0^p)$, with respect to the multi-index
$\mathbf{n}$ defined in \eqref{multi:index}. That is, for each $j=1,\ldots,p$, we have
\[
Q_{n}(\lambda) z_0^{j}(\lambda)-Q_{n-j} (\lambda)={\mathcal{O}}\Big(\frac{1}{\lambda^{n_{j}+1}}\Big),\qquad \lambda\rightarrow\infty,
\]
where $n_{j}$ is the $j$-th component of $\mathbf{n}$.
\item[iv)] For each $j=1,\ldots,p$
\[\lim_{n\to \infty} \frac{Q_{n-j}(\lambda)}{Q_n(\lambda)} = z_0^j(\lambda)\]
uniformly on each compact subset of $\overline{\mathbb{C}} \setminus \Gamma_1$.
\end{itemize}
\end{teo}

In particular, the last assertion of the theorem implies that the zeros of the polynomials $Q_n$ have $\Gamma_1$ as the set of its accumulation points. We would like to know more about this set and the functionals $L_j$ defined above, at least under additional assumptions on the coefficients $a_k$. When $a_0=\ldots= a_{p-1} = 0, a_p \neq 0$, the sets $\Gamma_k$ are starlike and the functionals $L_j$ admit integral representation in terms of (complex) rotationally  invariant measures  living on $\Gamma_1$ (see \cite{AKS}, \cite{LopezGarcia} and \cite{DL}). Moreover, the measures representing the functionals $L_j$ have the structure of a Nikishin system (see below for the definition and \cite{kn:Nik}). The zeros of the polynomials $Q_n$ are the eigenvalues of the $n$-th principal section of $A$. In \cite{kn:DK}, the asymptotic distribution of these eigenvalues and other points called generalized eigenvalues is obtained by solving a vector equilibrium problem on the sets $\Gamma_k, k=1,\ldots,p$. The structure of the vector equilibrium problem is the same as that appearing in finding the asympototic behavior of type II multiple orthogonal polynomials with respect to a Nikishin system of measures. These are only a few indicators which allow to predict that there is an underlying Nikishin structure in the study of  \eqref{rec} and even more general recurrence relations (see \cite{DLL}). In the next result we point out a new situation where this is so.

Set
\[ r(z) := a(1/z).\]
The algebraic equation $r(z) = \lambda$ is equivalent to the characteristic equation associated with \eqref{rec} and determines the same Riemann surface $\mathcal{R}$ as $a(z) = \lambda$. If $z(\lambda)$ denotes the inverse of $a(z)$, then $1/z(\lambda)$ is the inverse of $r(z)$. Notice that
\[ a'(z) = 0 \quad \Leftrightarrow \quad r'(1/z) = 0.\]
Therefore, the set of critical points of $a$ coincides with the reciprocals of the critical points of $r$.
We have
\begin{equation} \label{r} \lambda = r(z) = z + a_0 + \cdots + \frac{a_p}{z^p} = \frac{z^{p+1} + a_0z^p+\cdots+a_p}{z^p}.
\end{equation}
In this case, $\lambda = \infty$ is a branch point of order $p-1$ and a simple pole of the inverse of $r$.  The finite branch points are given by those $\lambda_0 = r(z_0)$ such that $r'(z_0)=0$ for some $z_0 \in {\mathbb{C}} \setminus \{0\}$.  The set of finite branch points of the algebraic equations $a(z) = \lambda, r(z) = \lambda$ is the same.

Define
\[ q(z) = z^{p+1} -a_1 z^{p-1} - 2a_2 z^{p-2} - \cdots - p a_p .
\]
Obviously
\[ r'(z) = \frac{q(z)}{z^{p+1}},
\]
and
\[ r'(z) =0\quad  \Leftrightarrow\quad q(z) =0.
\]

Now, we restrict our attention to the case when all the zeros of $q$ are real and simple. In particular, $a_j \in {\mathbb{R}}, j=1,\ldots,p$. We suppose also that $a_0 \in {\mathbb{R}}$.
We will assume that $p >1$; otherwise, if $p=1$ then $a_1> 0$ so that $q$ has real zeros and we fall in the case of the standard Chebyshev polynomials of the first kind for which all is known.

In the sequel $\{x_1,\ldots,x_{p+1}\}$ denotes the set of zeros of $q$. Notice that
\[ x_1+\cdots + x_{p+1} =0.
\]
Consequently, there has to be at least one positive and one negative zero of $q$. We will assume that $q$ has exactly $p$  zeros of one sign and $1$ of the other. One case reduces to the other substituting in the algebraic equation $z$ by $-z$.  When there are $p$ negative critical points of $r$ we  enumerate them as follows
\[ x_1 < \cdots < x_p < 0 < x_{p+1}.
\]
When there are $p$ positive critical points it is more convenient to index them so that
\[ x_{p+1} < 0 < x_p < \cdots < x_{1}.
\]
In either cases
\[ \{\lambda_1,\ldots,\lambda_{p+1}\}, \qquad \lambda_k = r(x_k), \qquad k=1,\ldots,p+1.
\]
denotes the collection of finite branch points of the algebraic equation $r(z) = \lambda$ (or, $a(z) = \lambda$).

Define the system of absolutely continuous measures $(\rho_1,\ldots,\rho_p)$ given by
\begin{equation}
 \label{generators}
\frac{d \rho_1(x)}{dx} = \frac{{\rm Im}(z_{0,-}(x))}{\pi}, \quad x \in \Gamma_1,\qquad \frac{d \rho_j(x)}{dx} = \frac{{\rm Im}(z_{j-1,-}(x))}{\pi(x-\lambda_j)}, \quad x \in \Gamma_j,\quad j=2,\ldots,p.
 \end{equation}
For each $j=1,\ldots,p-1$ the product $\langle \beta_j , \beta_{j+1} \rangle$  of two measures $\beta_j,\beta_{j+1}$ supported on $\Gamma_j$ and $\Gamma_{j+1}$, respectively, is defined by
\begin{equation}
\label{prod}
d\langle \beta_j , \beta_{j+1} \rangle (x) = (x- \lambda_{j+1})\widehat{\beta}_{j+1}(x) d\beta_j(x), \qquad x \in \Gamma_j,
\end{equation}
where
\[\widehat{\beta}_{j+1}(\lambda):= \int \frac{d \beta_{j+1} (t)}{\lambda - t}\]
denotes the Cauchy transform of   $\beta_{j+1}$.
Finally, set
\[ \sigma_1 = \rho_1, \quad \sigma_2 = \langle \rho_1, \rho_2\rangle , \quad \sigma_3 = \langle \rho_1, \langle \rho_2,\rho_3\rangle\rangle , \ldots  \qquad \sigma_p = \langle \rho_1,\langle\rho_2\ldots, \rho_{p} \rangle \rangle.
\]
We will show that in the present situation all these measures are supported on subintervals of the real line and have constant sign. We say that  $(\sigma_1,\ldots,\sigma_p)$ is a (generalized) Nikishin system generated by $(\rho_1,\ldots,\rho_p)$ and denote this by writing $(\sigma_1,\ldots,\sigma_p) = {\mathcal{N}}(\rho_1,\ldots,\rho_p)$.

In standard Nikishin systems the product in \eqref{prod} is defined without the factor $x-\lambda_{j+1}$. It has been introduced here by necessity to compensate for the denominator $x-\lambda_{j}$ in \eqref{generators} which in turn is needed in order that $\widehat{\rho}_j(\lambda), j=2,\ldots,p,$ be well defined (note that $\rho_j$ is not finite). That the product of measures defined above is consistent in our case is discussed during the proof of

\begin{teo} \label{teo:2} Assume that the polynomial $q$ has $p$ zeros of one sign and $1$ of the other. With the indexing adopted above, we have:
\begin{itemize}
\item[i)] When the polynomial $q$ has $p$ negative zeros $\Gamma_1 = [\lambda_1,\lambda_{p+1}]$  and $\Gamma_1 = [\lambda_{p+1},\lambda_{1}]$ if it has $p$ positive zeros. For $k=2,\ldots,p$, $\Gamma_k = [\lambda_k,+\infty]$ if $x_k$ is a local maximum of $r$ and $\Gamma_k = [-\infty,\lambda_{k}]$ when $x_k$ is a local minimum of $r$.
\item[ii)] There exist a system of (signed) measures $(\mu_1,\ldots,\mu_m)$ supported on the interval $\Gamma_1$ absolutely continuous with respect to the Lebesgue measure such that
\begin{equation} \label{layer1}
g_k(\lambda) = \widehat{\mu}_k(\lambda) = \int \frac{d \mu_k(x)}{\lambda -x}, \qquad k=1,\ldots,p.
\end{equation}
\item[iii)] There exist real constants $c_{j,k}, j=1,\ldots,p, k=1,\ldots,j, c_{j,j} = 1,$ such that
\begin{equation} \label{sj}  \mu_j = \sum_{k=1}^{j} c_{j,k} \sigma_k, \qquad j=1,\ldots,p,
\end{equation}
where $(\sigma_1,\ldots,\sigma_p) = {\mathcal{N}}(\rho_1,\ldots,\rho_p)$ is the Nikishin system defined above.
\item[iv)] For each $j=1,\ldots,p$
\[ \int x^{k} Q_n(x) d\sigma_j(x) = 0, \qquad k= 0,\ldots,n_j -1.\]
The $n$ zeros of $Q_n$ are simple and lie in the interior of $\Gamma_1$.
\end{itemize}

\end{teo}

\section{Proof of Theorem \ref{teo:1}}

We divide the proof of Theorem \ref{teo:1} into several parts.

{\bf Proof of Theorem  \ref{teo:1} i).} Define the operations of division and multiplication of vectors as follows
\[ \frac{(1,1,\ldots,1)}{(y_1,y_2,\ldots,y_p)} := \left(\frac{1}{y_p},\frac{y_1}{y_p},\ldots, \frac{y_{p-1}}{y_p}\right),
\]
\[ (x_1,x_2,\ldots,x_p)(y_1,y_2,\ldots,y_p) := (x_1y_1,x_2y_2,\ldots,x_py_p).
\]
These operations define what is called the generalized Jacobi-Perron algorithm.  Applied iteratively, this algorithm allows to express a vector as a vector continued fraction. In \cite[Theorem 2, Remark 2]{Kal2} it was proved that the generalized Jacobi-Perron algorithm, applied to $(g_1,\ldots,g_p)$ outside the spectrum of $A$, gives the vector continued fraction
\[ (g_1,\ldots,g_p) = \frac{(1,1,\ldots,1)}{|(0, \ldots,0,\lambda - a_0)} + \cdots +
\frac{(1,1,\ldots,1)}{|(-a_{p-1},-a_{p-2},\ldots,-a_1,\lambda - a_0)} +
\]
\[\frac{(-a_p,1,\ldots,1)}{|(-a_{p-1},-a_{p-2},\ldots,-a_1,\lambda - a_0)} + \frac{(-a_p,1,\ldots,1)|}{|(-a_{p-1},-a_{p-2},\ldots,-a_1,\lambda - a_0)} + \cdots.
\]
Notice that the first $p$ floors of this continued fraction differ, but from the floor $p+1$ on, which we call the tail, they are identical. For brevity, we refrain from writing the variable $\lambda$ except where it is absolutely necessary to avoid confusion.

Let us  obtain an expression of the tail $(h_1,\ldots,h_p)$ of this continued fraction in terms of $h_1$. Obviously,
\[ (h_1,\ldots,h_p) = \frac{(-a_p,1,\ldots,1)}{|(-a_{p-1},-a_{p-2},\ldots,-a_1,\lambda - a_0) + (h_1,\ldots,h_p)}.
\]
Due to the way in which division and multiplication of vectors is defined we get that
\begin{equation} \label{aches} (h_1, \ldots,h_p) = \left(\frac{-a_p}{h_p + \lambda-a_0},\frac{h_1-a_{p-1}}{h_p + \lambda-a_0},\ldots, \frac{h_{p-1}-a_1}{h_p + \lambda-a_0}  \right).
\end{equation}

Identifying the first components we obtain
\begin{equation}  \label{eq:phi1} -\frac{h_1}{a_p} = \frac{1}{h_p + \lambda-a_0}.
\end{equation}
Equating the rest of the components, respectively, and using \eqref{eq:phi1} it follows that
\begin{equation} \label{eq:rec} h_k =  \frac{h_{k-1} - a_{p-k+1}}{h_p + \lambda -a_0} = -\frac{h_1}{a_p}(h_{k-1} - a_{p-k+1}),\qquad k=2,\ldots,p.
\end{equation}

These formulas allow to express all the components   in terms of $h_1$. Indeed, when $k=2$ we get
\[ h_2 = -a_p \left(\frac{-h_1}{a_p}\right)^2
- a_{p-1}\frac{-h_1}{a_p}.
\]
Substituting this in the second equality of \eqref{eq:rec}, for $k=3$ it follows that
\[ h_3 =  \frac{-h_1}{a_p}\left(-a_p \left(\frac{-h_1}{a_p}\right)^2
 -a_{p-1}\frac{-h_1}{a_p} - a_{p-2}\right) = -a_p\left(\frac{-h_1}{a_p}\right)^3   - a_{p-1}\left(\frac{-h_1}{a_p}\right)^2  - a_{p-2}\frac{-h_1}{a_p}.
\]
Applying this recursively we obtain
\begin{equation} \label{eq:phip} h_k = -a_p\left(\frac{-h_1}{a_p}\right)^{k} - a_{p-1}\left(\frac{-h_1}{a_p}\right)^{k-1}- \cdots -a_{p-k+1} \left(\frac{-h_1}{a_p}\right),\qquad k=2,\ldots,p.
\end{equation}
Combining \eqref{eq:phi1} with \eqref{eq:phip} for $k=p$, after some trivial manipulations we see that
\begin{equation} \label{ecalreves} a_p \left(\frac{-h_1}{a_p}\right)^{p+1} + a_{p-1} \left(\frac{-h_1}{a_p}\right)^{p} + \cdots + a_1 \left(\frac{-h_1}{a_p}\right)^2 + (a_0 - \lambda)\left(\frac{-h_1}{a_p}\right) +1 \equiv 0.
\end{equation}

This means that $  {-h_1}/a_p$ is one of the solutions of the algebraic equation \eqref{alg2}.
Furthermore, using the definition of product and division of vectors, \eqref{eq:phi1}, and the first equality in \eqref{eq:rec}, we can undo the first $p$ floors of the continued fraction getting
\[ (g_1,\ldots,g_p) = \frac{(1,1,\ldots,1)}{|(0, \ldots,0,\lambda - a_0)} + \cdots +
\frac{(1,1,\ldots,1)}{|(-a_{p-1},-a_{p-2},\ldots,-a_1,\lambda - a_0)} + (h_1,\ldots,h_p) =
\]
\[   \frac{(1,1,\ldots,1)}{|(0, \ldots,0,\lambda - a_0)} + \cdots +
\frac{(1,1,\ldots,1)}{|(0,-a_{p-2},\ldots,-a_1,\lambda - a_0)} + \left(-h_1/a_p,h_2,\ldots,h_p\right) =
\]
\[   \frac{(1,1,\ldots,1)}{|(0, \ldots,0,\lambda - a_0)} + \cdots +
\frac{(1,1,\ldots,1)}{|(0,0,-a_{p-3},\ldots,-a_1,\lambda - a_0)} + \left(-h_1/a_p,(-h_1/a_p)^2,h_3,\ldots,h_p\right) =
\]
\[ \cdots = \left(-h_1/a_p,(-h_1/a_p)^2, \ldots,(-h_1/a_p)^p\right).
\]
Since $-h_1/a_p$ solves the algebraic equation the equality of the first components renders that $g_1 = -h_1/a_p = z_0$ because $g_1(\infty) =0$.
It is worth noting that this equality used in \eqref{eq:phip} implies that
\begin{equation} \label{eq:phip2} h_k = -a_p z_0^{k} - a_{p-1}z_0^{k-1}- \cdots -a_{p-k+1} z_0,\qquad k=1,\ldots,p,
\end{equation}
in a neighborhood of $\infty$. This completes the proof.
\hfill $\Box$

\medskip

{\bf Proof of Theorem  \ref{teo:1} ii)-iv).} Consider the difference equation
\begin{equation}\label{recqnpn}
\lambda y_{n}=y_{n+1}+a_{0}\, y_{n}+a_{1}\, y_{n-1}+\cdots+a_{p}\, y_{n-p},\qquad n\geq 0,
\end{equation}
with initial conditions
\begin{equation}\label{initcond}
\begin{array}{cccccccc}
n & = & -p & -p+1 & -p+2 & \ldots & -1 & 0\\ \hline p_{n}^{(1)} & = & 1 & 0 & 0 & \ldots & 0 & 0 \\ p_{n}^{(2)} & = & 0
& 1 & 0 & \ldots & 0 & 0 \\ p_{n}^{(3)} & = & 0 & 0 & 1 & \ldots & 0 & 0 \\
 & \vdots & & & & \vdots & &\\
p_{n}^{(p)} & = & 0 & 0 & 0 & \ldots & 1 & 0 \\ q_{n} & = & 0 & 0 & 0 & \ldots & 0 & 1
\end{array}
\end{equation}
Observe that for $n\geq 0$, $q_{n}$ is a
polynomial of degree $n$ and $p_{n}^{(j)}$ is a polynomial of degree $n-1$ for all $j\in\{1,\ldots,p\}$. In fact, $q_n(\lambda) = Q_n(\lambda), n\geq 0,$ since both polynomials satisfy the same difference equation with the same initial conditions.

Define
\begin{equation} \label{fis}
\phi_{j}(z):=g_{1}(z)+g_{2}(z)+\cdots+g_{j}(z),\qquad 1\leq j\leq p.
\end{equation}
Observe that
\begin{equation}
\label{series}
\phi_{j}(\lambda)=\sum_{n=0}^{\infty}\frac{\mu_{j,n}}{\lambda^{n+1}}, \qquad  \mu_{j,n} = L_{j}(\lambda^n),
\end{equation}
for all $\lambda\in {\mathbb{C}}$ sufficiently large.

In \cite[Theorem 1.1, Lemma 2]{kn:Kal} it was proved that for each $j=1,\ldots,p$
\[L_j(\lambda^{\nu} q_n) = 0, \qquad j=0,\ldots,n_j -1, \]
and
\begin{equation}
\label{HP}
q_{n}(\lambda) \phi_{j}(\lambda)-p_{n}^{(j)}(\lambda)=\mathcal{O}\Big(\frac{1}{\lambda^{n_{j}+1}}\Big),\qquad z\rightarrow\infty,
\end{equation}
where $n_{j}$ is the $j$-th component of $\mathbf{n}$. Since $q_n = Q_n$, ii) follows. (Incidentally, \cite[Theorem 1.1]{kn:Kal} also contains the  vector continued expansion of $(\phi_1,\ldots,\phi_p)$.)

Let $A_{n}$ be the truncation of the matrix $A$ to the first $n$ rows and columns ($n\geq0$). We define
\begin{align}
D_{n}(\lambda) & :=\det (\lambda I_{n}-A_{n}),\label{defDn}\\ D_{n}^{(j)}(\lambda) & :=\det (\lambda
I_{n-j}-\left(A_{n}\right)^{[1,\ldots,j;1,\ldots,j]}),\qquad 1\leq j\leq p, \label{defnDnj}
\end{align}
where $I_{k}$ denotes the identity matrix of size $k$, and $M^{[1,\ldots,j;1,\ldots,j]}$ denotes the submatrix of $M$
obtained after deleting the first $j$ rows and columns. Set
\[
D_{0}(\lambda)\equiv D_{1}^{(1)}(\lambda)\equiv D_{2}^{(2)}(\lambda)\equiv \cdots \equiv D_{p}^{(p)}(z)\equiv 1,\\
\]
and take $D_{n}(\lambda)\equiv 0$, $n<0$, and $D_{n}^{(j)}(\lambda)\equiv 0$, $n<j$.

Expanding $D_n$ by its last row, we have
$$
D_{n} = Q_n = q_n,\qquad n\geq 0.
$$
c:3In \cite[Lemma 2.1]{DLL} it was proved that the following relation holds for every $j\in\{1,\ldots,p\}$:
\begin{equation}\label{relpnDn}
p_{n}^{(j)}=D_{n}^{(1)}+D_{n}^{(2)}+\cdots+D_{n}^{(j)},\qquad n\geq 0.
\end{equation}
From  \eqref{fis}, \eqref{HP}, \eqref{relpnDn},  and Theorem  \ref{teo:1} i), it readily follows that
\[
D_{n}(\lambda) z_0^{j}(\lambda)-D_{n}^{(j)}(\lambda)=\mathcal{O}\Big(\frac{1}{\lambda^{n_{j}+1}}\Big),\qquad z\rightarrow\infty, \qquad j=1,\ldots,p.
\]
That is, $(D_{n}^{(1)}/D_n,\ldots,D_{n}^{(p)}/D_n)$ is a type II Hermite-Pad\'e approximation of $(z_0,\ldots,z_0^p)$ with respect to  $\bf n$. To complete the proof of iii) it remains to point out that since the diagonals of the matrix $A$ are constant we have  $D_n^{(j)} = Q_{n-j}$.

Finally, using a formula of Widom \cite[Theorem 2.8]{BG}, in \cite[Proposition 5.1]{kn:DK} it was shown that
\[ D_n(\lambda) = Q_n(\lambda) = \frac{-1}{a_pz_0^{n+1}(\lambda)}\prod_{j=1}^p(z_j(\lambda) - z_0(\lambda))^{-1}(1 + {\mathcal{O}}(\exp(-c_K n)))\qquad n\to \infty,
\]
uniformly on every compact subset $K$ of ${\mathbf{C}} \setminus \Gamma_1$, where $c_K$ is a positive constant which depends on $K$. Taking convenient ratios of these polynomials the previous formula immediately gives iv). \hfill $\Box$

\section{Proof of Theorem \ref{teo:2} i)}\label{sec:3}

In the rest of the paper, we restrict our attention to the case when the polynomial $q$ has only $p$ zeros of a given sign (positive or negative) and the remaining one has the opposite sign. As mentioned in the introduction, without loss of generality we can assume that $p$ zeros are negative and one positive.

\begin{figure}[h]
\vskip-.2truein
\subfigure[]{
\includegraphics[width=0.32\textwidth]{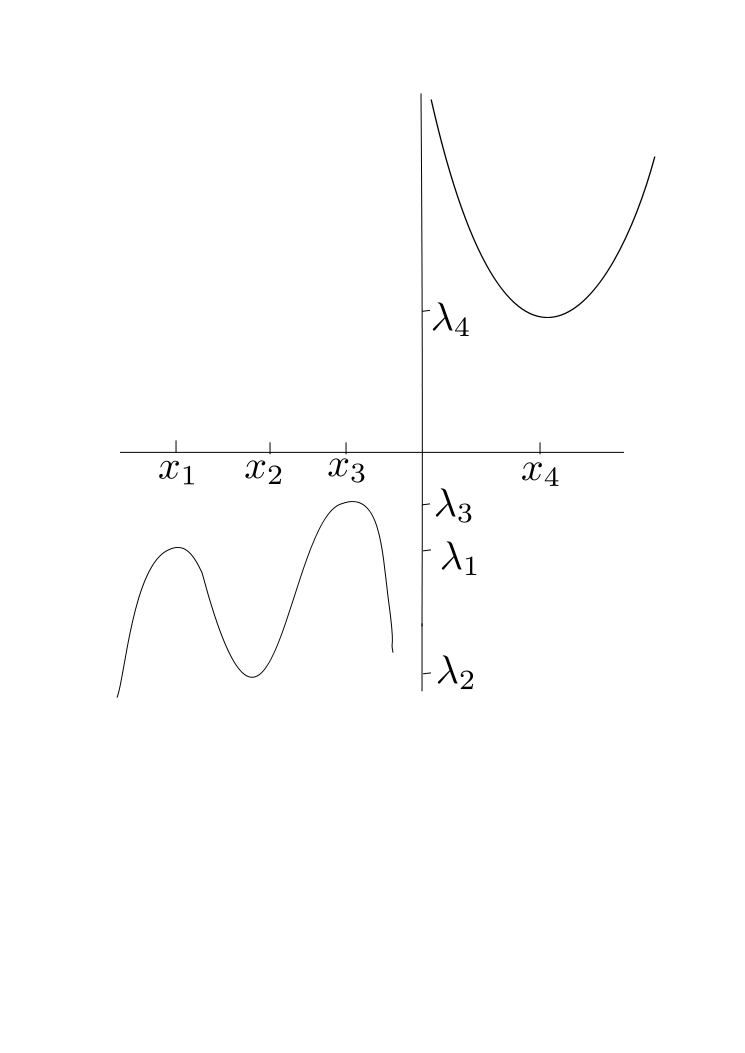}%\caption{}
}
\subfigure[]{
\includegraphics[width=0.32\textwidth]{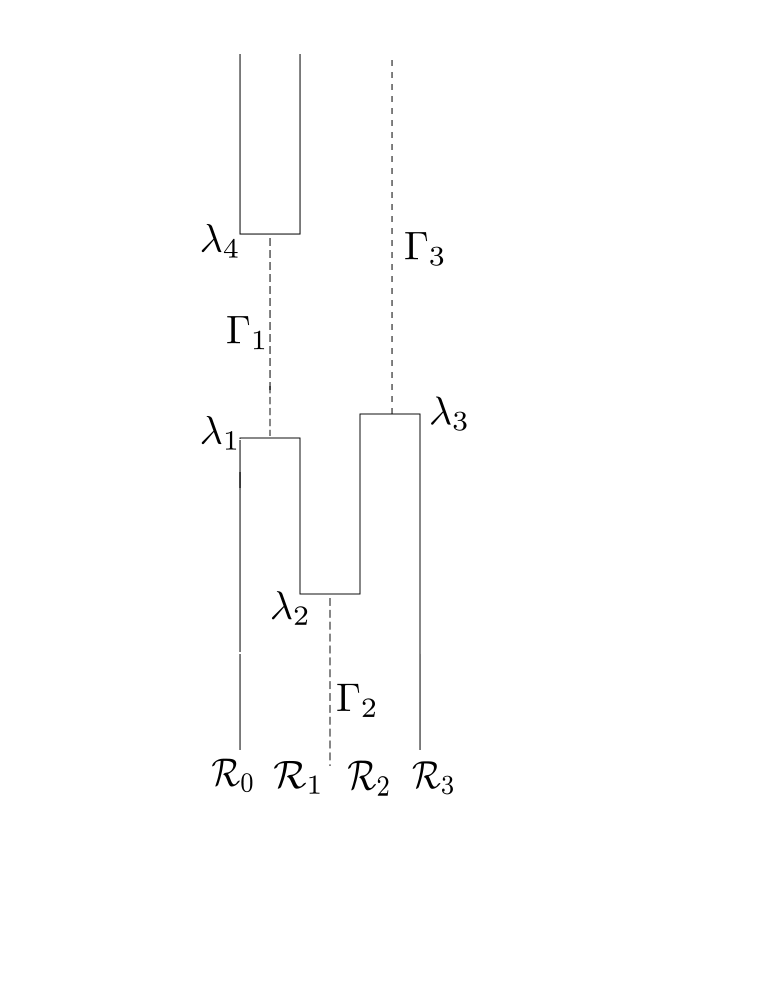}%\caption{}
}
\subfigure[]{
\includegraphics[width=0.32\textwidth]{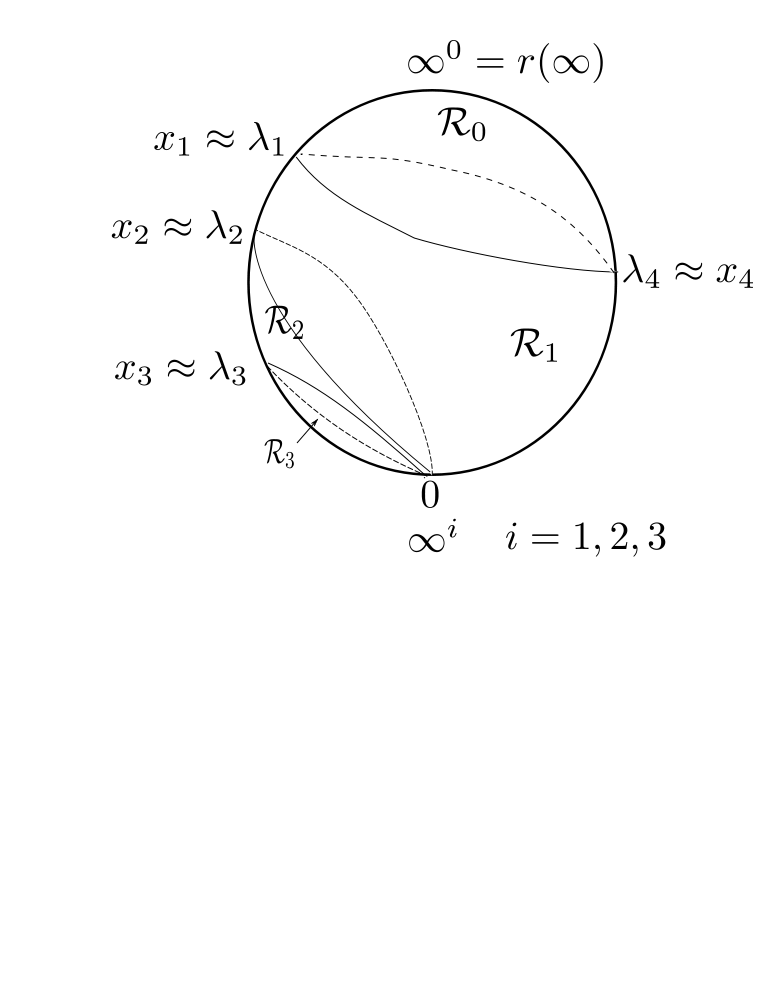}%\caption{}
}
%\end{subfigure}
\caption{(a) $p=3$, example with critical points $x_i$ and branch points $\lambda_i$;
(b) Placement of the  cuts $\Gamma$; (c) The cut sphere}
\end{figure}

Notice that  $1/z_0(\lambda)$ is holomorphic in a neighborhood of $\infty$ and has a simple pole at $\infty$. In the upper and lower half planes there are no branch points; therefore, $z_0(\lambda)$ can be extended analytically to the complement of a segment contained in the real line whose extreme points are branch points. On $(-\infty,x_1)$ and $(x_{p+1},+\infty)$ the rational function $r(z)$ is strictly increasing and invertible in a neighborhood of each point. Obviously, $\lambda_1 = r(x_1) < \lambda_{p+1} = r(x_{p+1})$ because  $r(-\infty,x_1] \cap r[x_{p+1},+\infty) = \emptyset$, since the algebraic equation defined by $r(z) = \lambda$ is irreducible. Define
\[ \tilde{\Gamma}_1 = [\lambda_1,\lambda_{p+1}],\qquad \tilde{\mathcal{R}}_0 = \overline{\mathbb{C}} \setminus \tilde{\Gamma}_1.
 \]
and let $\tilde{z}_0(\lambda)$ be the analytic extension of ${z}_0(\lambda)$ to $\tilde{\mathcal{R}}_0$.
Notice that $\lambda_1$ and $\lambda_{p+1}$ are first order branch points so $\tilde{z}_0$ returns to its initial value over any closed cycle contained in $\tilde{\mathcal{R}}_0$. Consequently $\tilde{z}_0$ is single valued. On $(x_1,x_2)$ and $(0,x_p)$ the function $r$ is strictly decreasing with a local minimum at $x_2$. By continuity, the segments $r(x_1,x_2) = (\lambda_2,\lambda_1)$ and $r(0,x_p) = (\lambda_p,+\infty)$ must belong to the next sheet attached to $\tilde{ \mathcal{R}}_0$ along $\tilde{\Gamma}_1 $. Set
\[ \tilde{\Gamma}_2 = [-\infty,\lambda_2],\qquad \tilde{\mathcal{R}}_1 = \overline{\mathbb{C}} \setminus (\tilde{\Gamma}_1 \cup \tilde{\Gamma}_2).
 \]
Paste $\tilde{\mathcal{R}}_0$ and $\tilde{\mathcal{R}}_1$ along $\tilde{\Gamma}_1$ in a crosswise manner as usual. Through $\tilde{\Gamma}_1$ we can extend $\tilde{z}_0$ analytically to all $\tilde{\mathcal{R}}_1$ and we denote the extension $\tilde{z}_1$. It is easy to see that $\tilde{z}_1$ returns to its initial value over any closed cycle contained in $\tilde{\mathcal{R}}_1$ so it is single valued. By continuity, the segment $r(x_2,x_3) = (\lambda_2,\lambda_3)$ must belong to the next sheet  and $x_3$ is a local maximum of $r$. Define
\[ \tilde{\Gamma}_3 = [\lambda_3,+\infty],\qquad \tilde{\mathcal{R}}_2 = \overline{\mathbb{C}} \setminus (\tilde{\Gamma}_2 \cup \tilde{\Gamma}_3).
 \]
Paste $\tilde{\mathcal{R}}_1$ and $\tilde{\mathcal{R}}_2$ along $\tilde{\Gamma}_2$ in a crosswise manner.  We can extend $\tilde{z}_1$ analytically to all $\tilde{\mathcal{R}}_2$ through $\tilde{\Gamma}_2$ and we denote the extension $\tilde{z}_2$.
We can repeat these arguments until we reach the interval $(x_p,0)$.

In general terms, if $x_k, 2\leq k \leq p$, is a local maximum of $r$, define
$\tilde{\Gamma}_k = [\lambda_k,+\infty)$, whereas if it is a local minimum, we take
$\tilde{\Gamma}_k = (-\infty,\lambda_k]$. Set
\[ \tilde{\mathcal{R}}_k = \left\{
\begin{array}{ll}
\overline{\mathbb{C}} \setminus \tilde{\Gamma}_1, & k = 0, \\
\overline{\mathbb{C}} \setminus (\tilde{\Gamma}_k \cup \tilde{\Gamma}_{k+1}), & 1 \leq k \leq p-1, \\
\overline{\mathbb{C}} \setminus \tilde{\Gamma}_p, & k = p, \\
\end{array}
\right.
\]
For each $k=1,\ldots,p$,  we can extend $\tilde{z}_{k-1}$ analytically to all $\tilde{\mathcal{R}}_k$ through $\tilde{\Gamma}_k$ and we denote the extension $\tilde{z}_k$. These extensions are all single valued and
\[ {\mathcal{R}} = \overline{ \cup_{k=0}^p\tilde{\mathcal{R}}_k}.
\]
From the uniqueness of analytic continuation we get that for each $\lambda \in {\mathbb{C}}$,  $\tilde{z}_0(\lambda), \ldots, \tilde{z}_p(\lambda) $ are a permutation of the roots  $ {z}_0(\lambda), \ldots, {z}_p(\lambda) $ of the algebraic equation $a(z) = \lambda$. If we prove that
\begin{equation}\label{ordering:roots:tilde}
|\tilde{ z}_{0}(\lambda)|\leq |\tilde{ z}_{1}(\lambda)|\leq \cdots\leq |\tilde z_{p}(\lambda)|,\qquad \lambda\in\mathbb{C},
\end{equation}
and
\begin{equation}\label{ordering:roots:strict} |\tilde{ z}_k(\lambda)|<|\tilde{ z}_{k+1}(\lambda)|, \qquad
\lambda\in{\mathbb{C}}\setminus\tilde{\Gamma}_{k+1},
\end{equation}
then we can conclude from \eqref{order} that $\tilde{ z}_k(\lambda)=z_k(\lambda), \lambda\in\tilde{\mathcal{R}}_k, k=0,\ldots,p,$ and $\Gamma_k = \tilde{\Gamma}_k, k=1,\ldots,p$. The rest of the proof is dedicated to verify \eqref{ordering:roots:tilde}--\eqref{ordering:roots:strict}. We do
this  following  the scheme employed in \cite{DLL} to prove a similar result.

For each $k=1,\ldots, p,$ define a measure $\tilde{s}_k$ on $\tilde{\Gamma}_k$ with density
\begin{equation}\label{rho:0}  d \tilde{s}_{k}(x):=\frac{1}{2\pi i}\left( \frac{\tilde
z_{k-1,+}'(x)}{\tilde z_{k-1,+}(x)}- \frac{\tilde z_{k-1,-}'(x)}{\tilde z_{k-1,-}(x)} \right) d  x = \frac{1}{\pi} \mbox{\rm Im}\left( \frac{\tilde
z_{k-1,+}'(x)}{\tilde z_{k-1,+}(x)} \right),\qquad x\in\tilde{\Gamma}_k,
\end{equation}
where the prime denotes the derivative with respect to $x$, and where the $+$ and $-$ subscripts stand
for the boundary values obtained from the upper or lower half of the complex plane, respectively, and $\mbox{\rm Im}(\cdot)$ stands for the imaginary part of a complex number $(\cdot)$. The second equality in \eqref{rho:0} comes from the fact that $z(\overline \lambda) = \overline{z(\lambda)}, \lambda \in {\mathcal{R}}$.

Notice that the density
\eqref{rho:0} is well-defined in the interior of $\tilde{\Gamma}_k$ ($z(\lambda)\neq 0, \lambda \in {\mathcal{R}}$). From the algebraic equation \eqref{alg2} it is easy to deduce that
\[ \tilde z_0(\lambda) =  \lambda^{-1}(1 + {\mathcal{O}}(\lambda^{-1})), \qquad \tilde z_k(\lambda) =  c_k\lambda^{1/p}(1 + {\mathcal{O}}(\lambda^{-1/p})), \qquad k=1,\ldots, p , \qquad \lambda \to \infty,\]
where $c_k$ is one of the roots of $c^p = 1/a_p$ taken in some order depending on the branch, and if $\tilde \lambda$ is a branch point of $\tilde z_k$
\[ \tilde z_k(\lambda) =  \tilde z_k(\tilde \lambda) +  {\mathcal{O}}((\lambda -\tilde \lambda)^{1/2})), \qquad k=0,\ldots, p , \qquad \lambda \to \tilde \lambda.\]
Therefore,
\begin{equation}
\label{eq:a}
 \frac{\tilde z_k'(\lambda)}{\tilde z_k(\lambda)} = \left\{
\begin{array}{lll}
-(1 + {\mathcal{O}}( {1}/{\lambda}))/{\lambda}, & \lambda \to \infty, & k=0, \\
(1 + {\mathcal{O}}( {1}/{\lambda^{1/p}}))/ (p\lambda), & \lambda \to \infty, & k=1,\ldots,p, \\
\end{array}
\right.\end{equation}
and
\begin{equation}
\label{eq:b}
\frac{\tilde z_k'(\lambda)}{\tilde z_k(\lambda)} =  {\mathcal{O}}((\lambda - \tilde \lambda)^{-1/2}), \qquad \lambda \to \tilde \lambda.\end{equation}
Consequently,
\begin{equation}
\label{eq:c}
 \frac{\tilde z_{k,+}'(\lambda)}{\tilde z_{k,+}(\lambda)} - \frac{\tilde z_{k,-}'(\lambda)}{\tilde z_{k,-}(\lambda)}= \left\{
\begin{array}{lll}
  {\mathcal{O}}( {1}/{\lambda^2}), & \lambda \to \infty, & k=0, \\
 {\mathcal{O}}({1}/{\lambda^{(p+1)/p}})), & \lambda \to \infty, & k=1,\ldots,p, \\
\end{array}
\right.\end{equation}
and
\begin{equation}
\label{eq:d}
\frac{\tilde z_{k,+}'(\lambda)}{\tilde z_{k,+}(\lambda)} - \frac{\tilde z_{k,-}'(\lambda)}{\tilde z_{k,-}(\lambda)} =  {\mathcal{O}}((\lambda - \tilde \lambda)^{-1/2}), \qquad \lambda \to \tilde \lambda.\end{equation}
In particular, \eqref{eq:c}-\eqref{eq:d} imply that the measures $\tilde s_1,\ldots,\tilde s_p$ are finite.

We claim that $\tilde{s}_k$ is a real-valued (possibly signed) measure on $\tilde{\Gamma}_k$ with total mass
\begin{equation}\label{rho:1} \tilde{s}_k(\tilde \Gamma_k) := \int_{\tilde{\Gamma}_k} d \tilde{s}_k(x)
= \frac{p-k+1}{p} ,\qquad k=1,\ldots,p.
\end{equation}
To prove equation~\eqref{rho:1}, using contour integration, we first derive the following relation for the Cauchy transforms
\begin{equation}\label{rho:2} \frac{\tilde z_k'(\lambda)}{\tilde z_k(\lambda)} =  -   \int_{\tilde{\Gamma}_k} \frac{ d \tilde s_{k} (x)}{x - \lambda} +  \int_{\tilde{\Gamma}_{k+1}} \frac{ d \tilde s_{k+1} (x)}{x - \lambda}, \qquad k=0,\ldots, p,
\end{equation}
where by convention $\tilde{\Gamma}_0 = \tilde{\Gamma}_{p+1} = \emptyset$ and $\tilde s_0 = \tilde s_{p+1}=0$.  We show in detail how this is done for $k=2,\ldots,p-1$. The cases $k=0,1,p$ have a slightly different geometric configuration but are in fact simpler. Notice that $\lambda_{k}, \lambda_{k+1}
$ are the branch points of $\tilde z_k$ on the boundary of $\tilde{\mathcal{R}}_k$ and to fix the geometry suppose that $\lambda_{k} < \lambda_{k+1}$.

Let $\gamma_R = \{\zeta: |\zeta| = R\}$ and $\gamma_{\varepsilon,j} = \{\zeta : |\zeta - \lambda_j| = \varepsilon\}$, where $R>0$ is sufficiently large and $\varepsilon>0$ is sufficiently small. $\gamma_R$ is oriented positively and $\gamma_{\varepsilon,k}, \gamma_{\varepsilon,k+1}$ negatively.  By Cauchy's integral formula, we have
\[ \frac{\tilde z_k'(\lambda)}{\tilde z_k(\lambda)} = \frac{1}{2\pi i} \int_{\gamma_R \cup \gamma_{\varepsilon,k} \cup \gamma_{\varepsilon,k+1}} \frac{\tilde z_k'(\zeta)}{\tilde z_k(\zeta)} \frac{  d\zeta}{\zeta - \lambda} -   \int_{-R}^{\lambda_{k} - \varepsilon} \frac{ d \tilde s_{k} (x)}{x - \lambda} +  \int^{R}_{\lambda_{k+1} + \varepsilon} \frac{ d \tilde s_{k+1} (x)}{x - \lambda},
\]
where $\lambda$ is any point lying inside the region limited by $\gamma_R, \gamma_{\varepsilon,k+1}, \gamma_{\varepsilon,k}, \tilde \Gamma_k$ and $\tilde \Gamma_{k+1}$. In the second integral it is used that $\tilde z_{k,\pm} (x) = \tilde z_{k-1,\mp} (x)$ on $\tilde \Gamma_k$.
From \eqref{eq:c}-\eqref{eq:d}, we have
\[  \lim_{\varepsilon \to 0, R \to +\infty} \left|\frac{1}{2\pi i} \int_{\gamma_R \cup \gamma_{\varepsilon,k} \cup \gamma_{\varepsilon,k+1}} \frac{\tilde z_k'(\zeta)}{\tilde z_k(\zeta)} \frac{  d\zeta}{\zeta - \lambda}\right| = 0,\]
uniformly on any compact subset of $\tilde{\mathcal{R}}_k$. Therefore, \eqref{rho:2} follows.

Multiplying \eqref{rho:2} by $\lambda$, taking the limit as $\lambda \to \infty$, using \eqref{eq:a}, and \eqref{eq:c}, we obtain
\begin{equation}\label{rho:3}   \tilde s_k (\tilde{ \Gamma}_{k}) - \tilde s_{k+1} (\tilde {\Gamma}_{k+1}) \left\{
\begin{array}{ll}
- 1,  & k=0, \\
\frac{1}{p}, & k=1,\ldots,p, \\
\end{array}
\right.
\end{equation}
(In passing the limit under the integral sign one can use Lebesgue's dominated convergence theorem due to \eqref{eq:a}.) Now, applying \eqref{rho:3} consecutively, starting with $k=0$ we obtain \eqref{rho:1}. (Recall that $\tilde {\Gamma}_0 = \tilde {\Gamma}_{p+1} = \emptyset$.)

In \cite[(4.16)]{DLL} a formula similar to \eqref{rho:1} is obtained which eventually leads to the proof of a statement similar to that of Theorem \ref{teo:1} i). In the situation considered in \cite{DLL} all the curves $\Gamma_k$ and intervals $\tilde \Gamma_k, k=1,\ldots,p$ (denoted there $\Delta_k$) are bounded and the coefficients of the recurrence relation \eqref{rec}, which depend on $n$, are periodic. Despite these differences, the rest of the arguments used in \cite[pp. 56-60]{DLL} to prove that $\Gamma_k =\tilde \Gamma_k, k=1,\ldots,p,$ can be transposed to our situation without any change. We limit ourselves to pointing out that aside from the measures $\tilde s_i,\ldots,\tilde s_p$ which live on the intervals $\tilde \Gamma_1, \ldots, \tilde \Gamma_p$ the authors also use the measures $s_1,\ldots, s_p$ with density
\begin{equation}
\label{meds}
 d s_k(x) := \frac{1}{2\pi i} \sum_{j=0}^{k-1} \left( \frac{
z_{j,+}'(x)}{ z_{j,+}(x)}- \frac  {z_{j,-}'(x)}{ z_{j,-}(x)} \right) d  x  ,\qquad x \in \Gamma_k, \qquad k=1,\ldots,p,
\end{equation}
which live on the curves $\Gamma_k$. (Notice that here we used the branches $z_0,\ldots,z_{p-1}$ instead of $\tilde z_0,\ldots,\tilde z_{p-1}$). The properties needed relative to the measures $s_1,\ldots, s_p$ and the curves $ \Gamma_1, \ldots,  \Gamma_p$ can be found in \cite[Theorem 2.2 (a), Proposition 3.2]{kn:DK}. We leave the details to the reader. \hfill $\Box$

\medskip
We wish to mention that in \cite[Theorem 2.3 (b)-(c)]{kn:DK} it is shown that the system of measures $(s_1,\ldots,s_p)$ is the solution of an interesting (Nikishin type) vector equilibrium problem on the system of curves  $(\Gamma_1,\ldots,\Gamma_p)$ (regardless of our assumptions on the critical points).

In conclusion, we have $\Gamma_k = \tilde \Gamma_k,  s_k = \tilde s_k, k=1,\ldots,p$ and $z_k(\lambda) = \tilde z_k(\lambda), \lambda \in {\mathcal{R}}_k= \tilde {\mathcal{R}}_k, k=0,\ldots,p$. These facts will be used in the sequel without further notice. Perhaps the equality $ s_k = \tilde s_k$ is not so obvious so we prove it.

Indeed the assertion is trivial for $k=1$. Since $z_0$ is holomorphic in ${\mathbf{C}} \setminus \Gamma_1$ and $\Gamma_1 \cap \Gamma_2 = \emptyset$
\[ \frac{
z_{0,+}'(x)}{ z_{0,+}(x)}- \frac  {z_{0,-}'(x)}{ z_{0,-}(x)} \equiv 0, \qquad x \in \Gamma_2, \]
and $ {s}_2 = \tilde s_2$. For $3 \leq k \leq p$, it is easy to see that $z_0\cdots z_{k-2}$ is holomorphic and different from zero in the simply connected region ${\mathbb{C}} \setminus \Gamma_{k-1}$ and $\Gamma_{k-1} \cap \Gamma_k = \emptyset$. Therefore, $\log (z_0\cdots z_{k-2})$ is holomorphic in a neighborhood of $\Gamma_k$ and
\[(\log (z_0\cdots z_{k-2}))_+'(x) - (\log (z_0\cdots z_{k-2}))_-'(x) \equiv 0 , \qquad \qquad x \in   \Gamma_k. \]
Consequently, $s_k = \tilde s_k$.

\section{Proof of Theorem \ref{teo:2} ii)-iv)}

The functions $g_k,k=1,\ldots,p$ are defined outside the spectrum of the resolvent operator $R_\lambda$, see \eqref{resolvf}. Under the present assumptions we know that they can be extended analytically to ${\mathbb{C}} \setminus \Gamma_1$. For simplicity, the extension will also be denoted $g_k$. Thus,
\[g_k(\lambda) =  {z}^k_0(\lambda), \qquad \lambda \in {\mathbb{C}} \setminus \Gamma_1. \]
Among other things, we wish to show that the functions $g_k$ are Cauchy transforms of measures supported on $\Gamma_1$. For this purpose, we will use the following lemma.

\begin{lemma} \label{laux} Let $g\in {\mathcal{H}}({\mathbb{C}} \setminus \Gamma)$ where $\Gamma$ is a semi-infinite interval contained in the real line. Assume that $g$ extends continuously to $\Gamma$ from the upper and lower half planes and the limiting values are denoted $g_+$ and $g_-$, respectively. Suppose that
\[ g(\lambda) = o(1), \qquad \lambda \to \infty,
\]
and
\[ g(\lambda) = {\mathcal{O}}( (\lambda - \lambda_0)^{\alpha}), \qquad \alpha > -1, \qquad \lambda \to \lambda_0,
\]
where $\lambda_0$ is the finite end point of $\Gamma$. Then
\[ g(\lambda) = \frac{1}{2\pi i} \int_{\Gamma} \frac{ (g_-(x)- g_+(x)) dx}{\lambda - x}.
\]
\end{lemma}

{\bf Proof.} Let $\gamma_R = \{\zeta: |\zeta| = R\}$ and $\gamma_{\varepsilon} = \{\zeta : |\zeta - \lambda_0| = \varepsilon\}$, where $R>0$ is sufficiently large, $\varepsilon>0$ is sufficiently small, and the curves are oriented in the positive and negative directions, respectively. To fix ideas, assume that $\Gamma = (-\infty,\lambda_0]$. Then, by Cauchy's integral formula, we have
\[ g(\lambda) = \frac{1}{2\pi i} \int_{\gamma_R \cup \gamma_{\varepsilon}} \frac{g(\zeta) d\zeta}{\zeta - \lambda} + \frac{1}{2\pi i} \int_{-R}^{\lambda_0 - \varepsilon} \frac{(g_+(x) - g_-(x))dx}{x - \lambda},
\]
where $\lambda$ is any point lying inside the region limited by $\gamma_R, \gamma_{\varepsilon},$ and $\Gamma$. Making $R \to \infty$, $\varepsilon \to 0$, and using the hypothesis the integral representation immediately follows. \hfill $\Box$

\medskip
Now we are ready to prove assertions ii)-iii) of Theorem \ref{teo:2}.   Since all the coefficients $a_k$ are real, the resolvent functions are symmetric with respect to the real line, thus
\[ g_k(\overline{\lambda}) = \overline{g_k(\lambda)}, \qquad k=1,\ldots,p.
\]
Additionally,
\[ g_k \in {\mathcal{H}}(\overline{\mathbb{C}} \setminus \Gamma_1), \qquad g_k(\lambda) = {\mathcal{O}}\left(\frac{1}{\lambda}\right), \qquad \lambda \to \infty\qquad k=1,\ldots,p.
\]
Therefore, using the Cauchy integral theorem
\[ g_k(\lambda) = \frac{1}{2\pi i} \int_{\gamma} \frac{g_k(\zeta)}{\lambda - \zeta} d\zeta,
\]
for all $\lambda$ exterior to $\gamma$, where $\gamma$ is a positively oriented closed Jordan curve that surrounds $\Gamma_1$. Shrinking $\gamma$ to $\Gamma_1$, we find that \eqref{layer1} takes place with
\[ \frac{d \mu_k(x)}{d x} := \frac{\left( g_{k,-}(x) -g_{k,+}(x)\right)}{2\pi i}  = \frac{\mbox{\rm Im}\left( g_{k,-}(x) \right)}{\pi}, \qquad x \in \Gamma_1, \qquad k=1,\ldots,p.
\]
Thus, $\supp(\mu_k) = \Gamma_1$ and $\mu_k$ is absolutely continuous with respect to the Lebesgue measure. We are done with ii).

When $k=1$,   we get that
\[ \frac{d \mu_1(x)}{d x} = \frac{1}{\pi} \mbox{\rm Im}   ({ {z}_{0,-}(x)}) ,\qquad x \in \Gamma_1.
\]
This function is continuous on $\Gamma_1$ and can never equal zero in the interior of $\Gamma_1$ because if that should happen then $ {z}_{0,-}(x) = {z}_{0,+}(x)$, but that is impossible (except at the end points) because $ {z}(\lambda)$ is one to one whereas $x_-$ and $x_+$ are distinct points on $\mathcal{R}$. So, $\mu_1$ has constant sign on $\Gamma_1$.

Notice that
$\lim_{\lambda \to \infty} \lambda  {z}_0(\lambda)  =1$, so $ {z}_0$ transforms the lower half plane into the upper half plane ($ {z}_0$ only takes real values on ${\mathbb{R}}\setminus \Gamma_1$). Therefore, in fact, $\mu_1$ is a positive measure.

To prove iii) we use induction.   Set
\[ g_k^{(1)} := g_k \qquad \mbox{and} \qquad
\mu_k^{(1)} := \mu_k, \qquad k=1,\ldots,p.
\] We have proved that $\mu_1^{(1)}$ is a measure with constant sign on $\Gamma_1$ and we take
\[\sigma_1 = \mu_1^{(1)}. \]

Whenever possible,  to reduce the notation, we refrain from writing down the independent variable, and indicate separately where the formulas take place. For $k=1,\ldots,p$, on $\Gamma_1$ we have
\[ g_{k,-}^{(1)} - g_{k,+}^{(1)} =    { {z}_{0,-}^k}   -  { {z}_{0,+}^k}   =
 \left(  { {z}_{0,-}}  -   { {z}_{0,+}} \right)   \sum_{k_1 + k_2=k-1}   { {z}_{0,-}^{ k_1}}    { {z}_{0,+}^{k_2}} .
\]
Notice that on $\Gamma_1$
\[ \sum_{k_1 + k_2=k-1}   { {z}_{0,-}^{ k_1}}    { {z}_{0,+}^{k_2}} = \sum_{k_1 + k_2=k-1}  { {z}_{0,-}^{ k_1}}    { {z}_{1,-}^{k_2}}=
\sum_{k_1 + k_2=k-1}  { {z}_{1,+}^{ k_1}}    { {z}_{0,+}^{k_2}}=
\sum_{k_1 + k_2=k-1}  { {z}_{1,+}^{ k_2}}    { {z}_{0,+}^{k_1}}.
\]
Consequently, this function extends to a holomorphic function in ${\mathbb{C}}\setminus \Gamma_2$ and
\begin{equation} \label{phi2} g_k^{(2)} := \sum_{k_1 + k_2=k-1}   {z}_{0}^{ k_1}   {z}_{1}^{k_2} \in {\mathcal{H}}({\mathbb{C}}\setminus \Gamma_2), \qquad k=2,\ldots,p.
\end{equation}
We have also shown that
\begin{equation} \label{s1} \frac{d \mu_k^{(1)}}{d x} = \frac{g_{k,-}^{(1)} - g_{k,+}^{(1)}}{2\pi i} =     \frac{g_k^{(2)}{\rm Im}({ {z}_{0,-}})}{\pi }, \qquad \mbox{on}\,\, \Gamma_1.
\end{equation}
In particular
\[ g_k^{(2)} = \frac{g_{k,-}^{(1)} - g_{k,+}^{(1)}}{g_{1,-}^{(1)} - g_{1,+}^{(1)}}, \qquad \mbox{on}\,\, \Gamma_1, \qquad k=2,\ldots,p.
\]

Let us obtain an integral representation of the functions $ {g_k^{(2)}}, k=2,\ldots,p$, with respect to a measure supported on $\Gamma_2$. We know that
\[  {z}_1(\lambda) = {\mathcal{O}}\left( {\lambda^{1/p}}\right), \qquad \lambda \to \infty,
\]
and
\[  {z}_1(\lambda) =  {z}_1(\lambda_2) +  {\mathcal{O}}((\lambda -\lambda_2)^{1/2}), \qquad \lambda \to \lambda_2.
\]
From the expression of $g_k^{(2)}$ in \eqref{phi2}, it follows that for $k=2,\ldots,p$
\[ g_k^{(2)}(\lambda) = {\mathcal{O}}\left( {\lambda^{(k-1)/p}}\right), \qquad \lambda \to \infty,
\]
and
\[ g_k^{(2)}(\lambda) = g_k^{(2)}(\lambda_2) +  {\mathcal{O}}\left((\lambda -\lambda_2)^{1/2}\right), \qquad \lambda \to \lambda_2.
\]
Thus
\[ \frac{g_k^{(2)}(\lambda) - g_k^{(2)}(\lambda_2)}{\lambda - \lambda_2} = {\mathcal{O}}\left( {\lambda^{(k-p-1)/p}}\right), \qquad \lambda \to \infty,
\]
and
\[ \frac{g_k^{(2)}(\lambda) - g_k^{(2)}(\lambda_2)}{\lambda - \lambda_2} = {\mathcal{O}}\left((\lambda -\lambda_2)^{-1/2}\right), \qquad \lambda \to \lambda_2.
\]
From Lemma \ref{laux}, we obtain that
\begin{equation} \label{int2} g_k^{(2)}(\lambda) = g_k^{(2)}(\lambda_2) + \frac{\lambda - \lambda_2}{2\pi i} \int_{\Gamma_2} \frac{\left(g_{k,-}^{(2)}(x)- g_{k,+}^{(2)}(x)\right) dx}{(x-\lambda_2)(\lambda- x)}.
\end{equation}
We denote by $\mu_k^{(2)}$ the absolutely continuous measure on $\Gamma_2$ given by
\[ \frac{d \mu_k^{(2)}(x)}{dx} := \frac{\mbox{\rm Im}\left(g_{k,-}^{(2)}(x)\right)}{\pi (x-\lambda_2)}, \qquad x \in \Gamma_2, \qquad  k=2,\ldots,p.
\]

Let us restrict our attention to $\mu_2^{(2)}$. In this case, from \eqref{phi2},
\[ g_{2,-}^{(2)}(x)- g_{2,+}^{(2)}(x) =    { {z}_{1,-}(x)}-  { {z}_{1,+}(x)} , \qquad x \in \Gamma_2,
\]
since $ {z}_{0,-}(x) =  {z}_{0,+}(x), x \in \Gamma_2,$ and
\begin{equation} \label{s2} \frac{d \mu_2^{(2)}(x)}{dx} = \frac{ \mbox{\rm Im}\left( {z}_{1,-} (x)\right)}{\pi (x-\lambda_2)}, \qquad x \in \Gamma_2,
\end{equation}
has constant sign on $\Gamma_2$. Indeed, if $\mbox{\rm Im}\left( {z}_{1,-} (x)\right) = 0$ at some interior point $x$ of $\Gamma_2$, then $ {z}_{1} (x_-) = {z}_{1} (x_+)$ which is impossible because $z_1$ is one-to-one. Since $\Gamma_1$ lies to one side of $\Gamma_2$, the integral
\[  \frac{\lambda - \lambda_2}{2\pi i} \int_{\Gamma_2} \frac{\left(g_{2,-}^{(2)}(x)- g_{2,+}^{(2)}(x)\right) dx}{(x-\lambda_2)(\lambda- x)}, \qquad \lambda \in \Gamma_1,
\]
also takes constant sign on $\Gamma_1$.

Putting together \eqref{s1}, \eqref{int2}, and \eqref{s2}, we obtain that
\[  \mu_2 = \mu_{2}^{(1)} = g_{2}^{(2)}(\lambda_2) \sigma_1 + \sigma_2,
\]
where $\sigma_2$, given by
\[ \frac{d\sigma_2(x)}{dx} =  (x - \lambda_2)  \int_{\Gamma_2} \frac{d\mu_2^{(2)}(t)}{ (x- t)} \frac{d\sigma_1(x)}{dx},
\]
has constant sign on $\Gamma_1$.

For each $2 \leq j \leq p$ and $j \leq k \leq p$, we define
\begin{equation} \label{fjk}
g_k^{(j)} :=  \sum_{k_1 + \cdots + k_j=k-j+1}   { {z}_0^{k_1}} { {z}_1^{k_2}}\cdots  { {z}_{j-1}^{k_j}},
\end{equation}
where $k_1,\ldots,k_j$ are non negative integers, and
\begin{equation} \label{sjk}
\frac{d \mu_k^{(j)}(x)}{dx} := \frac{\mbox{\rm Im}\left(g_{k,-}^{(j)}(x)\right)}{\pi (x-\lambda_j)}, \qquad x \in \Gamma_j, \qquad  k=j,\ldots,p.
\end{equation}
Then, we will show that $g_k^{(j)} \in {\mathcal{H}}\left({\mathbb{C}}\setminus \Gamma_j\right), j\leq k \leq p,$ and
\begin{equation} \label{fijk}
g_k^{(j)}(\lambda) = g_k^{(j)}(\lambda_j) + (\lambda - \lambda_j)\int_{\Gamma_j} \frac{d \mu_k^{(j)}(x)}{\lambda - x}, \qquad \lambda \in  {\mathbb{C}}\setminus \Gamma_j.
\end{equation}
In particular,
\[ g_j^{(j)}(\lambda) =    { {z}_0(\lambda)} + \cdots +  { {z}_{j-1}(\lambda)} , \qquad \lambda \in  {\mathbb{C}}\setminus \Gamma_j,
\]
and
\begin{equation} \label{jj} \frac{d \mu_j^{(j)}(x)}{dx} = \frac{  \mbox{\rm Im}\left( {z}_{j-1,-} (x)\right)}{\pi (x-\lambda_j)}, \qquad x \in \Gamma_j, \qquad j=2,\ldots,p,
\end{equation}
is a  measure with constant sign on $\Gamma_j$. For $j+1 \leq k \leq p$, we also have
\begin{equation} \label{sjk1}     {g_{k,-}^{(j)} - g_{k,+}^{(j)}}  = \left( { {z}_{j-1,-}} -   { {z}_{j-1,+}}\right)  {g_k^{(j+1)}} , \qquad \mbox{on}\,\, \Gamma_j,
\end{equation}
and
\begin{equation} \label{coc}  {g_k^{(j+1)}} = \frac{g_{k,-}^{(j)} - g_{k,+}^{(j)}}{g_{j,-}^{(j)} - g_{j,+}^{(j)}} , \qquad \mbox{on}\,\, \Gamma_j, \qquad k=j+1,\ldots,p.
\end{equation}

Let us show that \eqref{fijk}-\eqref{coc} are satisfied. Using induction we start proving that $g_k^{(j)} \in {\mathcal{H}}\left({\mathbb{C}}\setminus \Gamma_j\right)$, $j=1,\ldots,p, j\leq k\leq p$. We proved previously that this is true when $j=1,2$. Let us assume that the statement holds up to some $j-1, 2 \leq j-1 \leq p-1,$ and $k, j-1 \leq k \leq p,$ and prove that it is also verified for $j$ and $j\leq k\leq p$. By definition (we take $g_{j-2}^{(j-1)} \equiv 1$)
\begin{equation} \label{gjk} g_k^{(j)} = \sum_{k_j =0}^{k-j+1}  { {z}_{j-1}^{k_j}} \sum_{k_1 +\cdots+k_{j-1} = k -k_j -j+1}  { {z}_0^{k_1}}\cdots { {z}_{j-2}^{k_{j-1}}} = \sum_{k_j =0}^{k-j+1}  { {z}_{j-1}^{k_j}} g_{k-k_j-1}^{(j-1)}.
\end{equation}
By induction $g_{k-k_j-1}^{(j-1)} \in {\mathcal{H}}({\mathbb{C}}\setminus \Gamma_{j-1})$ and, obviously, $  { {z}_{j-1}^{k_j}} \in {\mathcal{H}}({\mathbb{C}}\setminus (\Gamma_{j-1}\cup \Gamma_{j}))$, so we only have to verify that $g_k^{(j)} \in {\mathcal{H}}(\Gamma_{j-1})$. This is done checking that $g_{k,-}^{(j)}(x) = g_{k,+}^{(j)}(x), x \in \Gamma_{j-1}$.

Let us rewrite $g_k^{(j)}$ as follows (again $g_{j-3}^{(j-2)} \equiv 1$)
\begin{equation} \label{gjk2} g_k^{(j)} = \sum_{k_j =0}^{k-j+1}  { {z}_{j-1}^{k_j}} \sum_{k_{j-1} =0}^{k-k_j - j+1}  { {z}_{j-2}^{k_{j-1}}} g_{k-k_j-k_{j-1}-2}^{(j-2)} = \sum_{\ell =0}^{k-j+1} g_{k-\ell-2}^{(j-2)}\sum_{k_{j-1}+k_{j} = \ell}  { {z}_{j-2}^{k_{j-1}}} { {z}_{j-1}^{k_j}}.
\end{equation}
By the induction hypothesis, we have $g_{k-k_j-k_{j-1}-2}^{(j-2)} \in {\mathcal{H}}({\mathbb{C}}\setminus \Gamma_{j-2} )$, in particular on $\Gamma_{j-1}$ its $\pm$ values are equal. So we must only worry about the $\pm$ values of the remaining part of the expression. We have
\[ \sum_{k_{j-1} + k_j =\ell}  { {z}_{j-2,-}^{k_{j-1}}}  { {z}_{j-1,-}^{k_j}}=
  \sum_{k_{j-1} + k_j =\ell}  { {z}_{j-1,+}^{k_{j-1}}}  { {z}_{j-2,+}^{k_j}} =
 \sum_{k_{j-1} + k_j =\ell}  { {z}_{j-2,+}^{k_{j-1}}}  { {z}_{j-1,+}^{k_j}},
\]
where the first equality comes from analytic extension through $\Gamma_{j-1}$ of $ {z}_{j-1}$ and $ {z}_{j-2}$, whereas the second equality is due to the symmetry of the indicated sum.  Thus the statement is true.

In order to derive \eqref{fijk}  it is sufficient to point out that
\[ \frac{g_{k}^{(j)}(\lambda)-g_{k}^{(j)}(\lambda_j)}{\lambda - \lambda_j} = \left\{
\begin{array}{ll}
{\mathcal{O}}\left(1/\lambda^{(k-j+1)/p}\right), &  \lambda \to \infty, \\
{\mathcal{O}}\left((\lambda - \lambda_j)^{-1/2}\right), &  \lambda \to \lambda_j,
\end{array}
\right.
\]
and use Lemma \ref{laux}. Formula \eqref{jj} follows directly from \eqref{sjk}, the definition of $g_j^{(j)}$, and the fact that $ { {z}_0 } + \cdots + { {z}_{j-2} } \in {\mathcal{H}}({\mathbb{C}}\setminus \Gamma_{j-1})$. Then, $\mu_j^{(j)}$ has constant sign on $\Gamma_j$ because  ${\rm Im}\left( { {z}_{j-1,-}}\right)  $ and $x - \lambda_j$ have fixed sign on that interval.

To get \eqref{sjk1} notice that from \eqref{gjk} we have on $\Gamma_j$
\[ g_{k,-}^{(j)} - g_{k,+}^{(j)}=   \sum_{k_j =0}^{k-j+1} \left( { {z}_{j-1,-}^{k_j}} -   { {z}_{j-1,+}^{k_j}}\right)g_{k-k_j-1}^{(j-1)} =
\]
\[\left( { {z}_{j-1,-}} -   { {z}_{j-1,+} }\right)
  \sum_{k_j =1}^{k-j+1} \sum_{\ell_1 + \ell_2 = k_j -1}  { {z}_{j-1,-}^{\ell_1}}  { {z}_{j-1,+}^{\ell_2}} g_{k-k_j-1}^{(j-1)}.
\]
Making the change of variables $k_j -1 = \ell$ and taking into consideration that $ {z}_{j-1,+} =  {z}_{j,-}$ on $\Gamma_j$, it follows that
\[ g_{k,-}^{(j)} - g_{k,+}^{(j)} = \left( { {z}_{j-1,-}} -   { {z}_{j-1,+} }\right)
  \sum_{\ell =0}^{k-j} g_{k-\ell-2}^{(j-1)} \sum_{\ell_1 + \ell_2 = \ell}  { {z}_{j-1,-}^{\ell_1}}  { {z}_{j,-}^{\ell_2}}.
\]
Since $g_{k}^{(j+1)} \in {\mathcal{H}}({\mathbb{C}} \setminus \Gamma_{j+1})$, according to \eqref{gjk2} (substitute $j$ by $j+1$) the last double sum coincides with $g_{k}^{(j+1)}$ and we are done. Finally, \eqref{coc} follows directly from \eqref{sjk1}.

Now, $   { {z}_{j,-}}  -   { {z}_{j,+}} $ has constant sign on $\Gamma_{j+1}$, thus
\[ (\lambda - \lambda_{j+1})\int_{\Gamma_{j+1}} \frac{d \mu_{j+1}^{(j+1)}(x)}{\lambda - x}
\]
has constant sign on the interval $\Gamma_j$ which is disjoint from $\Gamma_{j+1}$. For each $j=1,\ldots,p-1$ define a measure $\langle \mu_{j}^{(j)} , \mu_{j+1}^{(j+1)}\rangle$ whose differential expression is
\[ d\langle \mu_{j}^{(j)}, \mu_{j+1}^{(j+1)}\rangle (x) = (x- \lambda_{j+1})\widehat{\mu}_{j+1}^{(j+1)}(x) d\mu_j^{(j)}(x), \qquad x \in \Gamma_j.
\]
Each one of these measures has constant sign on its support.
Now, set
\[ \sigma_1 = \mu_{1}^{(1)}, \quad \sigma_2 = \langle \mu_{1}^{(1)}, \mu_{2}^{(2)}\rangle , \quad \sigma_3 = \langle \mu_{1}^{(1)}, \langle \mu_{2}^{(2)},\mu_{3}^{(3)} \rangle \rangle, \cdots \qquad \sigma_p = \langle \mu_{1}^{(1)},\ldots, \mu_{p}^{(p)} \rangle.
\]
Since $\Gamma_j \cap \Gamma_{j+1} = \emptyset, j=1,\ldots,p-1$, from the definition  it readily follows that
the measures $\sigma_1,\ldots,\sigma_p$ have constant sign on $\Gamma_1$. Obviously, $\mu_1 = \sigma_1$, $\mu_2 = g_2^{(2)}(\lambda_2)\sigma_1 + \sigma_2$, and the proof of \eqref{sj} follows directly employing \eqref{fijk} and \eqref{sjk1} several times. The constant $c_{j,k}, k=1,\ldots,j-1$ can be expressed in terms of the values $g_{k}^{(k)}(\lambda_k),k=1,\ldots,j-1$. Notice that the measures $\mu_j^{(j)}, j=1,\ldots,p,$ are the measures $\rho_j$ defined in \eqref{generators}. We have obtained iii).

From what has been proved (see \eqref{resolvf}, \eqref{L}, and \eqref{layer1}), we have for all $n \in {\mathbb{Z}}_+$
\[ L_j(x^n) = \sum_{k=1}^j \int x^n d \mu_k(x).
\]
Therefore, from   Theorem \ref{teo:1} ii), it follows that
\[ \int x^{\nu}Q_n(x) d\mu_j(x) =0 ,\qquad \nu=0,\ldots,n_j-1, \qquad j=1,\ldots,p,
\]
Due to \eqref{sj} and the decreasing structure of the values of the components of $(n_1,\ldots,n_p)$, we obtain
\begin{equation} \label{orto} \int x^{\nu}Q_n(x) d\sigma_j(x) =0 ,\qquad \nu=0,\ldots,n_j-1, \qquad j=1,\ldots,p.
\end{equation}
This settles assertion iv) and we conclude the proof. \hfill $\Box$

\section{Location of zeros and convergence of the Hermite-Pad\'e approximants}

Let us define the measures (see \eqref{prod})
\[ \rho_{j,j}= \rho_{j} , \qquad \rho_{j,k} = \langle \rho_{j} ,\ldots,\rho_{k} \rangle, \qquad 1 \leq j < k \leq p,
\]
where the measures $\rho_j$ are given by \eqref{generators}.
Let ${\bf n} = (n_1,\ldots,n_p), n_1 \geq n_2\geq \cdots \geq n_p,$ be a multi-index (not necessarily staircase). For each $j=1,\ldots,p$, consider the linear form
\begin{equation} \label{lf1} {\mathcal{L}}_{{\bf n},j}(\lambda) = \ell_j(\lambda)  + \sum_{k=j+1}^p \ell_k(\lambda)  (\lambda-\lambda_{j+1}) \widehat{\rho}_{j+1,k}(\lambda),\qquad  \deg \ell_k \leq n_k -1,
\end{equation}
(the sum is empty when $j=p$).

\begin{lemma}
\label{lem:5.1} With the notations set above, for each $j=1,\ldots,p$ the linear form $\mathcal{L}_{{\bf n},j}$ has in $\mathbb{C} \setminus \Gamma_{j+1} (\Gamma_{p+1} = \emptyset)$ at most $n_j +\cdots+n_p -1$ zeros or it is the null function.
\end{lemma}

{\bf Proof.} Obviously, the statement is true if $j=p$, so we can assume that $1\leq j <p$. Fix $j$. Assume that $\mathcal{L}_{{\bf n},j}$  is non null and  has at least $n_j +\cdots +n_p$ zeros in ${\mathbb{C}}\setminus \Gamma_{j+1}$. Let us prove that the reduced linear form $\mathcal{L}_{{\bf n},j+1}$
has at least $n_{j+1}+\cdots+n_p$ zeros in ${\mathbb{C}}\setminus \Gamma_{j+2}$.

The zeros of  ${\mathcal{L}}_{{\bf n},j}$ are symmetric with respect to the real line. Therefore, there exists a polynomial with real coefficients $w_{\bf n}$ of degree $\geq n_j +\cdots +n_p$ and zeros in ${\mathbb{C}}\setminus \Gamma_{j+1}$ such that
\[ \frac{{\mathcal{L}}_{{\bf n},j}}{w_{\bf n}} \in {\mathcal{H}}({\mathbb{C}}\setminus \Gamma_{j+1}).\]
Thus, for each $\nu = 0,\ldots,n_{j+1}+\cdots+n_p -1$
\[ \frac{\lambda^{\nu}{\mathcal{L}}_{{\bf n},j}}{w_{\bf n}} \in {\mathcal{H}}({\mathbb{C}}\setminus \Gamma_{j+1}),
\]
where $\deg \left(\lambda^{\nu}\ell_j(\lambda)\right) \leq n_j+\cdots+n_p-2$  and $\deg \left(\lambda^{\nu}(\lambda-\lambda_{j+1}) \ell_k(\lambda)\right) \leq n_j+\cdots+n_p-1.$

Let $\Gamma$ be a Jordan curve contained in ${\mathbb{C}}\setminus \Gamma_{j+1}$ with winding number $1$ that surrounds all the zeros of $w_{\bf n}$. Using Cauchy's theorem and Cauchy's integral formula, we get
\[ 0 = \int_{\Gamma} \frac{\lambda^{\nu}{\mathcal{L}}_{{\bf n},j}(\lambda)}{w_{\bf n}(\lambda)} d\lambda = \sum_{k=j+1}^p \int_{\Gamma}\frac{\lambda^{\nu} \ell_k(\lambda)  (\lambda-\lambda_{j+1})}{w_{\bf n}(\lambda)} \int \frac{d\rho_{j+1,k}(x)}{\lambda- x} d\lambda =
\]
\[ \sum_{k=j+1}^p \int \int_{\Gamma}\frac{\lambda^{\nu} \ell_k(\lambda)  (\lambda-\lambda_{j+1})d\lambda}{w_{\bf n}(\lambda)(\lambda- x)}  d\rho_{j+1,k}(x) = \int x^{\nu} \sum_{k=j+1}^p  \ell_k(x)  (x-\lambda_{j+1})   \frac{d\rho_{j+1,k}(x)}{w_{\bf n}(x)} =
\]
\[\int x^{\nu} {\mathcal{L}}_{{\bf n},j+1}(x)  \frac{(x-\lambda_{j+1})d\rho_{j+1}(x)}{w_{\bf n}(x)}.
\]
This implies that ${\mathcal{L}}_{{\bf n},j+1}$ has at least $n_{j+1}+\cdots+n_p$ sign changes on $\Gamma_{j+1}$ since the measure $ {(x-\lambda_{j+1})d\rho_{j+1}(x)}/{w_{\bf n}(x)}$ has constant sign on $\Gamma_{j+1}$. If ${\mathcal{L}}_{{\bf n},j+1}\equiv 0$ then $\ell_{k}\equiv 0,k=j+1,\ldots,p$ and consequently ${\mathcal{L}}_{{\bf n},j} \equiv 0$, otherwise it could not have had $n_j+\cdots+n_p$ zeros in $\mathbb{C} \setminus \Gamma_j$. Continuing the process, we conclude that ${\mathcal{L}}_{{\bf n},p} = \ell_p$ has at least $n_p$ zeros in $\mathbb{C}$ which is only possible if $\ell_p \equiv 0$ and going backward we obtain that ${\mathcal{L}}_{{\bf n},j}\equiv 0$, against our initial assumption. So ${\mathcal{L}}_{{\bf n},j}$ has at most $n_j+\cdots+n_p$ zeros in $\mathbb{C} \setminus \Gamma_{j+1}$. \hfill $\Box$

In the sequel, we return to multi-indices of the form \eqref{multi:index} for which we write ${\mathcal{L}}_{{\bf n},j}= {\mathcal{L}}_{{n},j}$.

\begin{teo} \label{zeros} For all $n \in {\mathbb{Z}}_+$ the zeros of $Q_n$ are simple, lie in the interior of $\Gamma_1$, and interlace the zeros of $Q_{n+1}$. Additionally,
\begin{equation} \label{coc2}  \limsup_{n\to \infty} \left\|z_0^j  - \frac{Q_{n-j}}{Q_{n} }\right\|_{\mathcal{K}}^{1/(n+n_j)} \leq \|\varphi\|_{\mathcal{K}}, \qquad j=1,\ldots,p.
\end{equation}
uniformly on compact subsets of ${\mathbb{C}}\setminus \Gamma_1$, where   $\varphi$ denotes the conformal representation of $\overline{\mathbb{C}}\setminus \Gamma_1$ onto the interior of the unit circle such that $\varphi(\infty) = 0$ and $\varphi^{\prime}(\infty) > 0$.
\end{teo}

{\bf Proof.} From \eqref{orto} it follows that for any system of polynomials $\ell_j, \deg \ell_j \leq n_j -1, j=1,\ldots,p$, with real coefficients
\begin{equation} \label{orto2} \int Q_n(x){\mathcal{L}}_{{n},1}(x) d\sigma_1(x) = 0.
\end{equation}
Then, \eqref{orto2} implies that $Q_n$ has exactly $n$ simple zeros in the interior of $\Gamma_1$.

Indeed, assume that $Q_{n}$ has at most $N \leq n-1$ sign changes in  the interior of $\Gamma_1$. Solving a homogeneous linear system of $n-1$ equations on the $n$ coefficients of the polynomials $\ell_k$ we can construct a non-trivial linear form with a simple zero at each of the $N$ points where $Q_n$ changes sign and a zero of order $n-1-N$ at one of the extreme points of $\Gamma_1$. From Lemma \ref{lem:5.1}, the corresponding linear form can have no other zero in the interior of $\Gamma_1$ except those imposed. Substituting in \eqref{orto2} this linear form we reach a contradiction.

Let $A, B$ be real constants such that $|A| + |B| \neq 0$. From \eqref{orto} it follows that for any system of polynomials $\ell_j, \deg \ell_j \leq n_j -1, j=1,\ldots,p$, with real coefficients
\begin{equation} \label{orto3} \int \left( AQ_n(x) + BQ_{n+1}(x)\right) {\mathcal{L}}_{{n},1}(x)  d\sigma_1(x) = 0.
\end{equation}
Using the property of the linear form ${\mathcal{L}}_{n,1}$ it follows that  $AQ_n(x) + BQ_{n+1}(x)$ has at least $n$ sign changes in the interior of $\Gamma_1$. Therefore, all the zeros of $AQ_n(x) + BQ_{n+1}(x)$ are simple and lie on the real line.

In order to prove that the zeros of $Q_{n+1}$ and $Q_n$ interlace let us show first that these polynomials do not have common zeros. To the contrary, suppose that there exists a point $x_{0} \in
\Gamma_1$ such that $Q_{n+1}(x_0) = Q_n(x_0) = 0$. As we have
seen, $x_{0}$ must be a simple zero of each one of these
functions. Therefore, $Q_n^{\prime}(x_{0}) \neq
0 \neq Q_{n+1}^{\prime}(x_{0})$. Thus, there must
exist real constants $A,B$ different from zero such that
\[  (AQ_n+BQ_{n+1})(x_{0}) = (AQ_n+BQ_{n+1})^{\prime}(x_{0}) = 0\,.
\]
This means that the function $AQ_n+BQ_{n+1} $ has a double zero at
$x_{0}$ which is impossible due to what was proved in the previous sentence.

Fix $y \in {\mathbb{R}} $ and set $G^y_n(z) = Q_{n+1}(z)Q_n(y)-Q_{n+1}(y)Q_n(z)$. That is, we take  $A = -Q_{n+1}(y)$ and $B = Q_n(y)$. Let $x_{\nu}, x_{\nu
+1}$ be two consecutive zeros of $Q_{n+1}$ in
$ \Gamma_1$ and let $y \in
(x_{\nu},x_{\nu+1})$. Then $|A| + |B| > 0.$

The function $G^y_n(z)$ is real  when
restricted to ${\mathbb{R}}$. We have $(G^y_n)^{\prime}(z) = Q_{n+1}^{\prime}(z)Q_n(y)- Q_{n+1}(y)Q_n^{\prime}(z)$.
Let us assume that $(G^{y_0}_{n})^{\prime}(y_0) = 0$ for some $y_0 \in (x_{\nu},x_{\nu-1})$. Since
$G^{y}_{n}(y) = 0$ for all $y \in (x_{\nu},x_{\nu+1})$
we obtain that $G^{y_0}_{n}(z)$ has a zero of order
$\geq 2$ (with respect to $z$) at $y_0$ which is impossible.  Consequently,
\[(G^{y}_{n})^{\prime}(y) =
Q_{n+1}^{\prime}(y)Q_n(y)- Q_{n+1}(y)Q_n^{\prime}(y)
\]
takes values with constant sign for all $y \in
(x_{\nu},x_{\nu+1})$. At the extreme points $x_{\nu}, x_{\nu+1}$
this function cannot be equal to zero because the polynomials
$Q_n, Q_{n+1}$ do not have common zeros.
By continuity, $(G^{y}_{n})^{\prime}$ preserves the same
sign on all $[x_{\nu},x_{\nu+1}]$ (and, consequently, on  $[x_1,x_{n+1}]$). Thus
\[ \mbox{sign}(G^{x_{\nu}}_{n})^{\prime}(x_{\nu}) =
\mbox{sign}(Q_{n+1}^{\prime}(x_{\nu})Q_n(x_{\nu})) =
\]
\[ \mbox{sign}(Q_{n+1}^{\prime}(x_{\nu+1})Q_n(x_{\nu +1})) =
\mbox{sign}(G^{x_{\nu + 1}}_{n})^{\prime}(x_{\nu +1})\,.
\]
Since
\[ \mbox{sign}Q_{n+1}^{\prime}(x_{\nu}) \neq \mbox{sign}Q_{n+1}^{\prime}(x_{\nu+1})\,,
\]
we obtain that
\[ \mbox{sign}Q_n(x_{\nu}) \neq \mbox{sign}Q_n(x_{\nu+1}).
\]
Hence, there must be an intermediate zero of $Q_n$ between $x_{\nu}$ and $x_{\nu+1}$.

Now that we know that the zeros of $Q_n$ lie in $\Gamma_1$, from Theorem \ref{teo:1} iii) and Theorem \ref{teo:2} i), for each $j=1,\ldots,p$
\[ z_0^{j}(z) - \frac{Q_{n-j} (z)}{Q_{n}(z)} = {\mathcal{O}}\left(\frac{1}{z^{n+n_j+1}}\right) \in {\mathcal{H}}(\overline{\mathbb{C}}\setminus \Gamma_1).
\]
Using Theorem \ref{teo:2} iv) (or the interlacing property of the zeros of the monic polynomials $Q_n$)  the family of functions $(Q_{n-1}/Q_n), n \in {\mathbb{N}}$, is uniformly bounded on each compact subset of $\overline{\mathbb{C}}\setminus \Gamma_1$. Fix a compact set ${\mathcal{K}} \subset {\mathbb{C}} \setminus \Gamma_1$ and fix $0 < \rho < 1$ sufficiently close to $1$ so that $\mathcal{K}$ lies in the unbounded connected component of the complement of $\Gamma_\rho = \{\lambda: |\varphi(\lambda)| = \rho\}$.
Then
\[\left(z_0^{j}(\lambda) - \frac{Q_{n-j} (\lambda)}{Q_{n}(\lambda)}\right)/\varphi^{n +n_j +1}(\lambda) \in {\mathcal{H}}({\mathbb{C}} \setminus \Gamma_1) \]
and
\[\left\|\left(z_0^{j}  - \frac{Q_{n-j} }{Q_{n} }\right)/\varphi^{n +n_j +1} \right\|_{\Gamma_\rho} \leq C(\Gamma_\rho)/\rho^{n+n_j+1}. \]
where $C(\Gamma_\rho)$ is a positive constant independent of $n$. Using the maximum principle it follows that
\[\left\| z_0^{j}  - \frac{Q_{n-j} }{Q_{n} } \right\|_{\mathcal{K}} \leq C(\Gamma_\rho)\left(\|\varphi\|_{\mathcal{K}}/\rho\right)^{n+n_j+1} . \]
Taking the $n+n_j+1$ root, $\limsup$ as $n\to \infty$, and letting $\rho \to 1$ we have \eqref{coc2}. \hfill $\Box$

\section{Second type functions and the intervals $\Gamma_j$}

Assume that $(\sigma_1,\ldots,\sigma_p) = {\mathcal{N}}(\rho_1,\ldots,\rho_p)$ is the (generalized) Nikishin system given in Theorem \ref{teo:2} and ${\bf n} = (n_1,\ldots,n_p)$ the multi-index \eqref{multi:index}. Let $Q_n$ be the $n$-th monic multiple orthogonal polynomial with respect to ${\mathcal{N}}(\rho_1,\ldots,\rho_p)$. Define recursively
\[ \Psi_{n,0} = Q_n,\quad \Psi_{n,1}(\lambda) = \int{\frac{\Psi_{n,0}(x) d\rho_1(x)}{\lambda-x}}, \quad \Psi_{n,j}(\lambda) = \int{\frac{\Psi_{n,j-1}(x) (x-\lambda_j)d\rho_j(x)}{\lambda-x}},\quad j=2,\ldots,p.
\]

\begin{lemma}\label{lem:6.1} For each $j=1,\ldots,p$
\begin{equation}
\label{order1}
\Psi_{n,j}(\lambda) = \mathcal{O}(1/\lambda^{n_j+2-j}), \qquad \lambda \to \infty,
\end{equation}
where the limit is taken along any curve non-tangential to $\Gamma_j$ at $\infty$.
\end{lemma}

{\bf Proof.} Let $q$ be an arbitrary polynomial of degree $ \leq n_1$. By Theorem \ref{teo:2}-iv)
\[
\int \frac{q(\lambda) - q(x)}{\lambda - x} Q_{n}(x) d\rho_1(x) = 0
\]
which means that
\[
 \Psi_{n,1}(\lambda) = \int \frac{Q_n(x) d\rho_1(x)}{\lambda -x}  = \int \frac{q(x) Q_n(x) d\rho_1(x)}{q(\lambda)(\lambda -x)} = \mathcal{O}\left(\frac{1}{\lambda^{n_1 +1}}\right), \qquad \lambda \to \infty.
\]
This is \eqref{order1} for $j=1$.

Moreover, using Theorem \ref{teo:2}-iv) for $k=2,\ldots,p$ and   $\nu = 0\ldots,n_k -1$, we have
\[ 0 = \int x^{\nu} Q_n(x)  d \langle \rho_1,\rho_{2,k}\rangle(x) = \int x^{\nu} Q_n(x) (x-\lambda_2) \widehat{\rho}_{2,k}(x) d  \rho_1 (x) = \]
\[ \int \int \frac{x^{\nu} Q_n(x) (x-\lambda_2)d\rho_1(x)}{x-t}  d \rho_{2,k}(t)  =  \int t^{\nu}  (t-\lambda_2)\int \frac{ Q_n(x) d\rho_1(x)}{x-t}  d \rho_{2,k}(t).\]
In the third equality, the definition of $\widehat{\rho}_{2,k}$ and Fubini's theorem are used. The conditions of Fubini's theorem hold because $\rho_1$ has  compact support and ${d \rho_{2,k}(t)}{dt} = \mathcal{O}(1/t^{1/p})$. In the fourth equality, we use that $\deg (t^{\nu}  (t-\lambda_2)) \leq n_k \leq n_1$ and that $Q_n$ is orthogonal to all polynomials of degree $\leq n_1 -1$ with respect to $\rho_1$. The second integral after the last equality is $\Psi_{n,1}$ so for each $k=2,\ldots,p$ and $\nu = 0,\ldots,n_k -1$
\[ 0 =  \int x^{\nu}   \Psi_{n,1}(x)  (x-\lambda_{2}) d \rho_{2,k}(x).
\]

Assume that for some fixed $j \in \{1,\ldots,p-1\}$ we have \eqref{order1} and for each $k= j+1,\ldots,p$ and $ \nu=0,\ldots,n_k -j$
\begin{equation}
\label{orto7}
0 =  \int x^{\nu}   \Psi_{n,j}(x)  (x-\lambda_{j+1}) d \rho_{j+1,k}(x).
\end{equation}
Obviously, using induction, we complete the proof if we show that then \eqref{order1} holds for $j+1$ and also prove that whenever $j+1 \leq p-1$ we have for $k=j+2,\ldots,p$ and $\nu = 0,\ldots,n_k- j-1$
\[ 0 =  \int x^{\nu}   \Psi_{n,j+1}(x) (t-\lambda_{j+2})  d \rho_{j+2,k}(t).
\]

In fact, \eqref{orto7} with $k=j+1$ implies that for any polynomial $q$ of degree $\leq n_{j+1} -j +1$
\[
\int \frac{q(\lambda) - q(x)}{\lambda - x}\Psi_{n,j}(x)  (x-\lambda_{j+1}) d \rho_{j+1}(x) = 0.
\]
That is
\[
 \Psi_{n,j+1}(\lambda) =   \int \frac{q(x) (x-\lambda_{j+1}) \Psi_{n,j}(x) d\rho_{j+1}(x)}{q(\lambda)(\lambda -x)}.
\]
Observe that
\[\frac{q(x) (x-\lambda_{j+1}) \Psi_{n,j}(x)}{(\lambda -x)} = \mathcal{O}(1/x),\qquad x\to \infty, \qquad x \in \Gamma_{j+1}
\]
(uniformly with respect to $\lambda$ non-tangential to $\Gamma_{j+1}$) and $d\rho_{j+1}/dx = \mathcal{O}(1/x^{(p-1)/p}), x\to\infty,$  so the last integral is well defined and, therefore, $\Psi_{n,j+1} = \mathcal{O}(1/\lambda^{n_{j+1}-(j+1)+2})$ as needed.

On the other hand, \eqref{orto7} with $k= j+2,\ldots,p$  (in case that $j+2 \leq p$) and  $\nu = 0\ldots,n_{k} -j -1$ gives
\[ 0 = \int x^{\nu} \Psi_{n,j}(x) (x-\lambda_{j+1}) d\langle \rho_{j+1},\rho_{j+2,k}\rangle(x)  =
\]
\[ \int x^{\nu} \Psi_{n,j}(x)(x-\lambda_{j+1})(x -\lambda_{j+2}) \widehat{\rho}_{j+2,k}(x) d  \rho_{j+1} (x) =  \]
\[ \int \int \frac{x^{\nu} (x-\lambda_{j+1})(x -\lambda_{j+2}) \Psi_{n,j}(x)  d\rho_{j+1}(x)}{x-t}  d \rho_{j+2,k}(t)  =
\]
\[ \int t^{\nu}  (t-\lambda_{j+2})\int \frac{ \Psi_{n,j}(x) (x-\lambda_{j+1}) d\rho_{j+1}(x)}{x-t}  d \rho_{j+2,k}(t).\]
as needed, since the second integral in the last line is $\Psi_{n,j+1}$. In the third equality, the use of Fubini's theorem is justified because
$x^{\nu} (x-\lambda_{j+1})(x -\lambda_{j+2}) \Psi_{n,j}(x) = \mathcal{O}(1/x), d\rho_{j+1}(x)/dx = \mathcal{O}(1/x^{(p-1)/p}), (x-t)^{-1} = \mathcal{O}(1/t)$ and $d \rho_{j+2,k}(t)/dt = \mathcal{O}(1/t^{1/p})$ when $x \to \infty, x \in \Gamma_{j+1}$ and $t\to \infty, t\in \Gamma_{j+2}$. Thus, we have proved the assertion and the proof is complete.     \hfill $\Box$

\medskip

We have the following Widom type formulas for the second type functions.

\begin{teo}
\label{Widom}
For each $\ell = 0,\ldots, p$, we have
\begin{equation}
\label{eq:Widom}
\Psi_{n,\ell}(\lambda) = \frac{-1}{a_p} \sum_{j=\ell}^{p} \frac{1}{\prod_{k = \ell, k \neq j}^p (z_k(\lambda) - z_j(\lambda))} \frac{1}{z_j^{n+1}(\lambda)}
\end{equation}
(the product is replaced by $1$ when $\ell = p$). Moreover,
\begin{equation}
\label{strong}
\lim_{n \to \infty} z_{\ell}^{n+1} (\lambda)\Psi_{n,\ell}(\lambda) = -\left(a_p \prod_{k = \ell+1}^p (z_k(\lambda) - z_\ell(\lambda))\right)^{-1}
\end{equation}
uniformly on each compact subset of $\mathbb{C} \setminus (\Gamma_\ell \cup \Gamma_{\ell+1}), (\Gamma_0 = \emptyset = \Gamma_{p+1})$.
\end{teo}

{\bf Proof.} Notice that formula \eqref{eq:Widom} can be computed for all $\lambda \in \mathbb{C} \setminus \mathcal{A}$ where $\mathcal{A}$ consists of the endpoints of the intervals $\Gamma_j, j=1,\ldots, p,$ which are the only points $\lambda$ where one may have $z_j(\lambda) = z_k(\lambda)$ with $j \neq k$. For $\ell = 0$, \eqref{eq:Widom} is Widom's formula (see \cite[Theorem 2.8]{BG}). For $\ell = 1,\ldots,p$, we will prove the formula by induction on $\ell$ and for $\ell =0$ we know it holds. Let us assume that \eqref{eq:Widom} holds true when $\ell$ is replaced by $\ell -1, \ell \in \{1,\ldots,p\}$ and we will prove that it is also satisfied for $\ell$.

The function $\Psi_{n,l}$ satisfies the following boundary value problem:
\begin{itemize}
\item[(a)] $\Psi_{n,\ell} \in \mathcal{H}(\mathbb{C} \setminus \Gamma_\ell)$.
\item[(b)] It verifies the jump relation
\[ \Psi_{n,\ell,-}(x) - \Psi_{n,\ell,+}(x) = 2i \Psi_{n,\ell-1}(x) \mbox{Im} (z_{\ell-1,-}(x)), \qquad x \in \Gamma_\ell .
\]
\item[(c)] $\Psi_{n,\ell}$ remains bounded near the endpoint(s) of $\Gamma_{\ell}.$
\item[(d)]
\[\Psi_{n,\ell}(\lambda) = \mathcal{O}(1/ \lambda^{n_\ell - \delta}), \qquad \lambda \to \infty,
\]
where $\delta$ is independent of $n$
\end{itemize}
Indeed, $({\rm a})$  and $(\rm b)$ follow from the definition of $\Psi_{n,\ell}$ and the Sokhotski-Plemelj formula. $(\rm c)$ follows from the behavior of $z_{\ell -1}$ near the finite endpoint(s) of $\Gamma_{\ell}$ and the analyticity of $\Psi_{n,\ell-1}$ in a neighborhood of such points. Finally, $(\rm d)$ is a consequence of \eqref{order1}.  Given $\Psi_{n,\ell -1}$, this boundary value problem has a unique solution.

To prove \eqref{eq:Widom} it suffices to show that its right hand side, which we denote $\widetilde{\Psi}_{n,\ell}$,  verifies a weaker version of the boundary value problem. Let $\mathcal{A} = \{\lambda_1,\ldots,\lambda_{p+1}\}$ be the set of finite branch points of $z(\lambda)$. Then,
\begin{itemize}
\item[(a')] $\widetilde \Psi_{n,\ell} \in \mathcal{H}(\mathbb{C} \setminus (\Gamma_\ell \cup \mathcal{A}))$.
\item[(b')] It verifies the jump relation
\[ \widetilde \Psi_{n,\ell,-}(x) - \widetilde \Psi_{n,\ell,+}(x) = 2i \widetilde \Psi_{n,\ell-1}(x) \mbox{Im} (z_{\ell-1}(x)), \qquad x \in \Gamma_\ell\setminus \mathcal{A}.
\]
\item[(c')] Near each $\hat \lambda \in \mathcal{A}$
\[\widetilde \Psi_{n,\ell}(\lambda) = \mathcal{O}((\lambda - \hat\lambda)^{-1/2}), \qquad \lambda \to \hat \lambda.\]
\item[(d')]
\[\widetilde \Psi_{n,\ell}(\lambda) = \mathcal{O}(1/ \lambda^{n_\ell}), \qquad \lambda \to \infty.
\]
\end{itemize}
In fact, assume that for some $\ell \in \{1,\ldots,p\}$ \eqref{eq:Widom} holds for $\ell -1$ and let us prove that it also holds for $\ell$. Suppose that (a')-(d') hold. From (a), (b), (a'), (b'), and the induction hypothesis it follows that
\[\Psi_{n,\ell} - \widetilde \Psi_{n,l} \in \mathcal{H}(\mathbb{C} \setminus \mathcal{A}).\]
Then, (c) and (c') imply that the difference is also holomorphic at the points in $\mathcal{A}$  so it is an entire function. Finally, from (d) and (d'), using Liouville's theorem, we have $\Psi_{n,\ell} - \widetilde{\Psi}_{n,\ell} \equiv 0$ as needed.

Let us check (a')-(d'). Obviously, $\widetilde{\Psi}_{n,\ell} \in \mathcal{H}(\mathbb{C} \setminus \mathbb{R})$. On the other hand, it is finite except possibly  where $z_j = z_k$ for $j\neq k$; that is, on points in $\mathcal A$, and $\widetilde{\Psi}_{n,\ell}$  is a symmetric function of $z_{\ell},\ldots,z_{p}$ which does not depend on $z_0,\ldots,z_{\ell -1}$. Therefore, it is analytic in $\mathbb{C} \setminus (\Gamma_\ell \setminus \mathcal{A})$.

Regarding the jump condition, we start our with the terms in
\[\frac{-1}{a_p} \sum_{j=\ell+1}^{p} \frac{1}{\prod_{k = \ell, k \neq j}^p (z_k(\lambda) - z_j(\lambda))} \frac{1}{z_j^{n+1}(\lambda)} =
\]
\[\frac{-1}{a_p} \sum_{j=\ell+1}^{p} \frac{z_{\ell-1}(\lambda) - z_j(\lambda)}{(z_{\ell-1}(\lambda) - z_j(\lambda))(z_\ell(\lambda) - z_j(\lambda))\prod_{k = \ell+1, k \neq j}^p (z_k(\lambda) - z_j(\lambda))} \frac{1}{z_j^{n+1}(\lambda)}. \]
Now,
\[\sum_{j=\ell+1}^{p} \frac{  z_j(\lambda)}{(z_{\ell-1}(\lambda) - z_j(\lambda))(z_\ell(\lambda) - z_j(\lambda))\prod_{k = \ell+1, k \neq j}^p (z_k(\lambda) - z_j(\lambda))} \frac{1}{z_j^{n+1}(\lambda)}\]
is symmetric in $z_{\ell+1},\ldots,z_p$ and $z_{\ell-1}, z_{\ell}$, and doesn't depend on $z_0,\ldots,z_{\ell-2}$ so its $-$ and $+$ boundary values on $\Gamma_\ell \setminus \mathcal{A}$ coincide. On the other hand,
\[\sum_{j=\ell+1}^{p} \frac{  1}{(z_{\ell-1}(\lambda) - z_j(\lambda))(z_\ell(\lambda) - z_j(\lambda))\prod_{k = \ell+1, k \neq j}^p (z_k(\lambda) - z_j(\lambda))} \frac{1}{z_j^{n+1}(\lambda)}\]
presents the same type of symmetry. In conclusion,
\[\left(\frac{-1}{a_p} \sum_{j=\ell+1}^{p} \frac{1}{\prod_{k = \ell, k \neq j}^p (z_k(x) - z_j(x))} \frac{1}{z_j^{n+1}(x)}\right)_- - \left(\frac{-1}{a_p} \sum_{j=\ell+1}^{p} \frac{1}{\prod_{k = \ell, k \neq j}^p (z_k(x) - z_j(x))} \frac{1}{z_j^{n+1}(x)}\right)_+ =
\]
\begin{equation}
\label{condb1}
 \frac{-1}{a_p} \sum_{j=\ell+1}^{p} \frac{2i \mbox{Im} (z_{\ell -1,-}(x))}{\prod_{k = \ell-1, k \neq j}^p (z_k(x) - z_j(x))} \frac{1}{z_j^{n+1}(x)}, \qquad x \in \Gamma_\ell \setminus \mathcal{A}.
\end{equation}

For the remaining term, we  have
\[ \frac{1}{\prod_{k = \ell +1}^p (z_k(\lambda) - z_\ell(\lambda))} \frac{1}{z_\ell^{n+1}(\lambda)} =  \frac{z_{\ell -1}(\lambda) - z_{\ell}(\lambda)}{(z_{\ell -1}(\lambda) - z_{\ell}(\lambda))\prod_{k = \ell +1}^p (z_k(\lambda) - z_\ell(\lambda))} \frac{1}{z_\ell^{n+1}(\lambda)},
\]
and  $z_{\ell -1,-}(x) - z_{\ell,-}(x) = z_{\ell -1,-}(x) - z_{\ell-1,+}(x) = 2i \mbox{Im}( z_{\ell -1,-}(x)) = - (z_{\ell -1,+}(x) - z_{\ell,+}(x))$. So
\[ \left(\frac{-1}{a_p} \frac{1}{\prod_{k = \ell +1}^p (z_k(x) - z_\ell(x))} \frac{1}{z_\ell^{n+1}(x)}\right)_- - \left(\frac{-1}{a_p} \frac{1}{\prod_{k = \ell +1}^p (z_k(x) - z_\ell(x))} \frac{1}{z_\ell^{n+1}(x)}\right)_+ =
\]
\[ \frac{-1}{a_p} \frac{2i \mbox{Im} z_{\ell -1}(x)}{(z_{\ell,+}(x) - z_{\ell-1,+}(x))\prod_{k = \ell +1}^p (z_k(x) - z_{\ell-1,+}(x))} \frac{1}{z_{\ell-1,+}^{n+1}(x)} +
\]
\[\frac{-1}{a_p}\frac{2i \mbox{Im} z_{\ell -1}(x)}{(z_{\ell -1,+}(x) - z_{\ell,+}(x))\prod_{k = \ell +1}^p (z_k(x) - z_{\ell,+}(x))} \frac{1}{z_{\ell,+}^{n+1}(x)} =
\]
\begin{equation}
\label{condb2}
\frac{-1}{a_p} \sum_{j=\ell-1}^{\ell} \frac{2i \mbox{Im} (z_{\ell -1}(x))}{\prod_{k = \ell-1, k \neq j}^p (z_k(x) - z_j(x))} \frac{1}{z_j^{n+1}(x)}, \qquad x \in \Gamma_\ell \setminus \mathcal{A}.
\end{equation}
Again the symmetry with respect to $z_{\ell +1},\ldots,z_p$ is used.
Having in mind \eqref{condb1}, \eqref{condb2}, and the induction hypothesis $(b')$ is obtained.

Finally, (d') is a consequence of the behavior of the functions $z_1,\ldots,z_p$ at $\infty$, while (c') follows from their behavior at the finite end points of the intervals $\Gamma_j$. Therefore \eqref{eq:Widom} is satisfied.

To prove \eqref{strong} it is sufficient to notice that on any compact subset $\mathcal{K} \subset \mathbb{C} \setminus (\Gamma_{\ell} \cup \Gamma_{\ell +1})$ we have
\[\sup_{\lambda \in \mathcal{K}} \left|\frac{z_{\ell}(\lambda)}{z_{j}(\lambda)}\right| = q_{\ell,j}(\mathcal{K}) < 1, \qquad j = \ell +1,\ldots,p, \]
and use \eqref{eq:Widom}. \hfill $\Box$

\medskip

For the case of standard Nikishin systems and generating measures with compact support, the following statements (and even slightly stronger ones) are known to hold (see \cite[Propositions 1, 2]{kn:Gonchar}).

\begin{teo}
\label{teo:3}
For each $j=1,\ldots,p-1$, the following orthogonality relations hold
\begin{equation}
\label{orto4}
\int x^{\nu} \Psi_{n,j}(x) (x-\lambda_{j+1}) d\rho_{j+1,k}(x) = 0, \qquad  k= j+1,\ldots,p,\qquad \nu = 0,\ldots,n_k -\delta_j -1,
\end{equation}
where $\delta_1 = 0$ and $\delta_j=1, j=2,\ldots,p-1$. Set $N_{n,j} = n_{j+1}+\cdots+n_p, j=0,\ldots,p-1$. For $j=0,1, \Psi_{n,j}$ has exactly $N_{n,j}$ zeros
in $\mathbb{C} \setminus \Gamma_j \,(\Gamma_0 = \emptyset)$ they are simple and lie in the interior of $\Gamma_{j+1}$. $\Psi_{n,p}$ has no zero in $\mathbb{C}\setminus \Gamma_p$. For $j=2,\ldots,p-1, \Psi_{n,j}$ has at most $N_{n,j} +j-1$ zeros in $\mathbb{C} \setminus \Gamma_j$ and at least $N_{n,j}-p+j$ sign changes on $\Gamma_{j+1}$. Let $Q_{n,j+1}$ be the monic polynomial whose roots are the zeros of $\Psi_{n,j}$ on $\mathbb{C} \setminus \Gamma_j$ (repeated according to multiplicity). We have
\begin{equation}
\label{orto8}
\int x^{\nu} Q_n(x) \frac{d\rho_{1}(x)}{Q_{n,2}(x)}= 0, \qquad \nu = 0,\ldots, N_{n,1} -1,
\end{equation}
and for $j=1,\ldots,p-1$
\begin{equation}
\label{orto5}
\int x^{\nu} \Psi_{n,j}(x) (x-\lambda_{j+1}) \frac{d\rho_{j+1}(x)}{Q_{n,j+2}(x)}= 0, \qquad \nu = 0,\ldots, n_{j+1} +\deg Q_{n,j+2} -2,
\end{equation}
where $Q_{n,p+1} \equiv 1$.
\end{teo}

{\bf Proof.} For $j=1$, formula \eqref{orto4} was verified in the initial induction step during the proof of Lemma \ref{lem:6.1}. For the other values of $j$, \eqref{orto4} is almost \eqref{orto7} except a higher number of orthogonality relations. However, due to Theorem \ref{Widom}, now we know that $\Psi_{n,j}(\lambda) = \mathcal{O}(1/\lambda^{n_j}), j=2,\ldots,p$, which allows to increase the number of orthogonality relations obtained in \eqref{orto7} up to the quantity specified in \eqref{orto4}. Indeed, notice that $x^{\nu}(x-\lambda_{j+1})(x - \lambda_{j+2})\Psi_{n,j}(x) = \mathcal{O}(1/x), x\to \infty, x\in \Gamma_{j+1},$ for all $\nu =0,\ldots,n_j-2$ (see last sentence in the proof of Lemma \ref{lem:6.1}).

Using \eqref{orto4} for $j=1$ and $k=2,\ldots,p$, we obtain
\begin{equation}
\label{orto6}
 \int \Psi_{n,1}(x) \mathcal{L}_{n,2}(x)   (x-\lambda_2)d\rho_2(x) = 0,
 \end{equation}
where $\mathcal{L}_{n,2}(x) =  \ell_2(x) + \sum_{k=3}^p \ell_k(x) (x-\lambda_3) \widehat{\rho}_{3,k}(x) , \deg \ell_k \leq n_k -1, k=2,\ldots, p$ and the sum is empty when $p=2$. Using Lemma \ref{lem:5.1} it readily follows that $\Psi_{n,1}$ has at least $N_{n,1}$ sign changes in $\Gamma_2$.

Let $\mathcal{Z}_{n,1}$ be the collection of all the zeros of $\Psi_{n,1}$ in $\mathbb{C} \setminus \Gamma_{1}$ (repeated according to multiplicity). Let us prove that $\mathcal{Z}_{n,1}$ has exactly $N_{n,1}$ points showing thus that all the zeros of $\Psi_{n,1}$ in $\mathbb{C} \setminus \Gamma_{1}$ are simple and coincide with the points where this function changes sign on $\Gamma_2$.

Assume that $\mathcal{Z}_{n,1}$ contains more than $N_{n,1}$ points. Then, there exists a monic polynomial $Q_{n,2}^*, \deg Q_{n,2}^* = N_{n,1}^* > N_{n,1},$ with real coefficients, whose zeros are points in $\mathcal{Z}_{n,1}$ such that
\[ \frac{\Psi_{n.1}}{Q_{n,2}^*} \in \mathcal{H}(\mathbb{C} \setminus \Gamma_1).\]
On the other hand, $\Psi_{n,1}(\lambda) = \mathcal{O}(1/\lambda^{n_1 +1}), \lambda \to \infty$, so
\[ \frac{\Psi_{n,1}(\lambda)}{Q_{n,2}^*(\lambda)} = \mathcal{O}\left(\frac{1}{\lambda^{n_1+ N_{n,1}^*+1}}\right), \qquad \lambda \to \infty.
\]
Let $\nu = 0,\ldots,n_1+ N_{n,1}^*-1 $. Choose a curve $\Gamma$ with winding number $1$ that surrounds $\Gamma_1$. Using Cauchy's and Fubini's theorems, we have
\[ 0= \frac{1}{2\pi i} \int_{\Gamma } \frac{\lambda^{\nu}}{Q_{n,2}^*(\lambda)} \int \frac{  Q_n(x) d\rho_1(x) }{ \lambda -x} d\lambda  =
\]
\[\int  Q_n(x) \frac{1}{2\pi i} \int_{\Gamma } \frac{\lambda^{\nu}}{Q_{n,2}^*(\lambda)}  \frac{d\lambda }{ \lambda -x}  d\rho_1(x) =  \int \frac{x^{\nu} Q_n(x)d\rho_1(x)}{Q_{n,2}^*(x)}. \]
This relation is clearly impossible  because  $N_{n,1}^* > N_{n,1}$ and $Q_n$ would be orthogonal to itself.
Therefore, $\mathcal{Z}_{n,1}$ contains exactly $N_{n,1}$ points. Define $Q_{n,2}$ as the monic polynomial of degree $N_{n,1}$ whose simple zeros are the points in $\mathcal{Z}_{n,1}$. Repeating the arguments above with $Q_{n,2}$ replacing $Q_{n,2}^*$ and taking $\nu=0,\ldots,N_{n,1}-1$, we obtain \eqref{orto8}.

Formulas \eqref{orto4} and \eqref{orto5} coincide when $j=p-1$.   Fix $j, 1 \leq j \leq p-2$. Assume that $\Psi_{n,j}$ has at most $N_{n,j} +j-1$ zeros in $\mathbb{C} \setminus \Gamma_j$ and at least $N_{n,j} -p +j$ sign changes on $\Gamma_{j+1}$. For $j=1$ the assumption is true. Let us prove that $\Psi_{n,j+1}$ has at most $N_{n,j+1} +j$ zeros in $\mathbb{C} \setminus \Gamma_{j+1}$ and at least $N_{n,j+1} -p +j +1$ sign changes on $\Gamma_{j+2}$.

Using \eqref{orto4} for the index $j+1$, we have
\[ \int \Psi_{n,j+1}(x) \mathcal{L}_{n,j+2}(x)  (x-\lambda_{j+2}) d\rho_{j+2}(x) = 0,
\]
where $\mathcal{L}_{n,j+2}(x) =  \ell_{j+2}(x) + \sum_{k=j+3}^p \ell_k(x) (x-\lambda_{j+3}) \widehat{\rho}_{j+3,k}(x) , \deg \ell_k \leq n_k -2, k=j+2,\ldots, p$ and the sum is empty if $j+3 > p$. From Lemma \ref{lem:5.1} it follows that $\Psi_{n,j+1}$ has at least $N_{n,j+1} -p +j +1$ sign changes on $\Gamma_{j+2}$.

Let $\mathcal{Z}_{n,j+1}$ be the collection of all the zeros of $\Psi_{n,j+1}$ in $\mathbb{C} \setminus \Gamma_{j+1}$ (counting multiplicities). Let us prove that $\mathcal{Z}_{n,j+1}$ contains at most $N_{n,j+1}+j$ points.

Assume that $\mathcal{Z}_{n,j+1}$ contains at least $N_{n,j+1} +j +1$ points. Then, there exists a polynomial $Q_{n,j+2}^*, \deg Q_{n,j+2}^* = N_{n,j+1}^* > N_{n,j+1}+j,$ with real coefficients, whose zeros are points in $\mathcal{Z}_{n,j+1}$ such that
\[ \frac{\Psi_{n.j+1}}{Q_{n,j+2}^*} \in \mathcal{H}(\mathbb{C} \setminus \Gamma_{j+1}).\]
On the other hand, $\Psi_{n,j+1}(\lambda) = \mathcal{O}(1/\lambda^{n_{j+1}}), \lambda \to \infty$, so
\[ \frac{\Psi_{n,j+1}(\lambda)}{Q_{n,j+2}^*(\lambda)} = \mathcal{O}\left(\frac{1}{\lambda^{n_{j+1}+ N_{n,j+1}^*}}\right), \qquad \lambda \to \infty.
\]
Let $\nu = 0,\ldots,n_{j+1}+ N_{n,j+1}^*-2 $. Choose a curve $\Gamma$ with winding number $1$ that surrounds all the zeros of $Q_{n,j+2}^*$.
Using Cauchy's  theorem and the definition of $\Psi_{n,j+1}$, we have
\[ 0= \frac{1}{2\pi i} \int_{\Gamma } \frac{\lambda^{\nu}}{Q_{n,j+2}^*(\lambda)} \int \frac{  \Psi_{n,j}(x) (x-\lambda_{j+1})d\rho_{j+1}(x) }{ \lambda -x} d\lambda.
\]
If $\nu = 0,\ldots,N_{n,j+1}^*-1$, from Fubini's theorem and Cauchy's integral formula, it follows that
\[0= \int  \Psi_{n,j}(x) (x-\lambda_{j+1}) \frac{1}{2\pi i} \int_{\Gamma } \frac{\lambda^{\nu}}{Q_{n,j+2}^*(\lambda)}  \frac{d\lambda }{ \lambda -x}  d\rho_{j+1}(x) =  \int \frac{ x^{\nu} \Psi_{n,j}(x) (x-\lambda_{j+1})d\rho_{j+1}(x) }{Q_{n,j+2}^*(x)}. \]
When $\nu = N_{n,j+1}^* -1 + \tilde{\nu}, \tilde{\nu} =1,\ldots, n_{j+1}-1$, we also employ \eqref{orto4} to get
\[0 = \frac{1}{2\pi i} \int_{\Gamma } \frac{\lambda^{\nu}}{Q_{n,j+2}^*(\lambda)} \int \frac{  \Psi_{n,j}(x) (x-\lambda_{j+1})d\rho_{j+1}(x) }{ \lambda -x} d\lambda = \]
\[\frac{1}{2\pi i} \int_{\Gamma } \frac{\lambda^{N_{n,{j+1}}^* -1}}{Q_{n,j+2}^*(\lambda)} \int \frac{ x^{\tilde{\nu}} \Psi_{n,j}(x) (x-\lambda_{j+1})d\rho_{j+1}(x) }{ \lambda -x} d\lambda =\]
\[\int x^{\tilde{\nu}} \Psi_{n,j}(x) (x-\lambda_{j+1}) \frac{1}{2\pi i} \int_{\Gamma } \frac{\lambda^{N_{n,{j+1}}^* -1}}{Q_{n,j+2}^*(\lambda)}  \frac{d\lambda }{ \lambda -x}  d\rho_{j+1}(x) =  \int \frac{ x^{\nu} \Psi_{n,j}(x) (x-\lambda_{j+1})d\rho_{j+1}(x) }{Q_{n,j+2}^*(x)}. \]

These relations imply that $\Psi_{n,j}$ has at least $n_{j+1}+ N_{n,j+1}^*-1 > N_{n,j} +j-1$ sign changes on $\Gamma_{j+1}$ which contradicts the induction hypothesis that $\Psi_{n,j}$ has at most $N_{n,j} +j-1$ zeros in all of $\mathbb{C} \setminus \Gamma_{j}$. Therefore, the induction is complete and the assertion about the zeros of $\Psi_{n,j}, j=1,\ldots,p-1$ follows. The absence of zeros of $\Psi_{n,p}$  in $\mathbb{C} \setminus \Gamma_p$ is a trivial consequence of formula \eqref{eq:Widom} for $\ell = p$.

In order to prove \eqref{orto5}, define $Q_{n,j+2}$ as the monic polynomial  whose roots are the points in $\mathcal{Z}_{n,j+1}$ and repeat the arguments employed above replacing $Q_{n,j+2}^*$ with $Q_{n,j+2}$.  We are done. \hfill $\Box$

\medskip

Let $T_n(z^{-k}(a(z) - \lambda))$ denote the Toeplitz matrix associated  with the symbol $z^{-k}(a(z) - \lambda), k=0,\ldots,p-1$.
More precisely, if $A_n$ is the $n$-th principal section of $A$ in \eqref{eq:A} and $I_n$ the identity matrix of order $n\times n$, then
$T_n(z^{-k}(a(z) - \lambda))$ is the matrix of order $n\times n$ obtained from $A_{n+k} - \lambda I_{n+k}$ after eliminating its first $k$ rows and last $k$ columns.

For each $k=0,\ldots,p-1$, set
\[ P_{n,k}(\lambda) := \mbox{det}\, T_n(z^{-k}(a(z) - \lambda)), \qquad \mbox{sp}_k\,T_n(a) := \{\lambda \in {\mathbb{C}}\,:\, P_{n,k}(\lambda)=0\}.
\]
In particular, $P_{n,0} = (-1)^n Q_n $. The set $\mbox{sp}_k\,T_n(a)$ is called the $k$-th generalized spectrum of $A_n$ and its points $k$-th generalized eigenvalues (so, $\mbox{sp}_0\,T_n(a)$ is the normal spectrum). These generalized spectra are intimately connected with the curves $\Gamma_{k+1}$.

Define
\[ \liminf_{n\to \infty} \mbox{sp}_k\,T_n(a)
\]
as the set of all points $\lambda \in {\mathbb{C}}$ such that there exists a sequence $\{\lambda_n\}_{n \in {\mathbb{N}}}, \lambda_n \in \mbox{sp}_k\,T_n(a),$ with $\lim_n \lambda_n = \lambda$, and
\[ \limsup_{n\to \infty} \mbox{sp}\,T_n(a)
\]
as the set  of all points $\lambda \in {\mathbb{C}}$ such that there exists a subsequence $\{\lambda_n\}_{n \in {\Lambda}}, \lambda_n \in \mbox{sp}_k\,T_n(a), \Lambda \subset {\mathbb{N}},$ and $\lim_{n\in \Lambda} \lambda_n = \lambda$. Also define
\[s_{n,k} = \frac{1}{n} \sum_{\lambda \in \mbox{sp}_k\,T_n(a)} \delta_{\lambda},\]
where in the sum each $\lambda$ is counted according to its multiplicity as a zero of $P_{n,k}$.

In \cite[Theorem 2.6]{kn:DK} it was proved that
\begin{equation}
\label{zeros2}
\liminf_{n\to \infty} \mbox{sp}_k\,T_n(a) = \limsup_{n\to \infty} \mbox{sp}_k\,T_n(a) = \Gamma_{k+1},\qquad k=0,\ldots,p-1,
\end{equation}
and
\begin{equation}
\label{weak}
\lim_{n\to \infty} \int \phi(\lambda) d s_{n,k}(\lambda) = \int \phi(\lambda) d s_{k+1}(\lambda),
\end{equation}
for every every bounded continuous function $\phi$ on $\mathbb{C}$, where the measures $s_{k+1}$ are those defined in \eqref{meds}. For $k=0$, this result is due to Schmidt and Spitzer \cite{SS}. We  point out that these results were obtained in \cite{kn:DK} with no restriction on the coefficients $a_k$ in \eqref{rec}. In the case we are considering, we have carried out numerous numerical experiments in which the generalized eigenvalues lie exactly on the corresponding segments $\Gamma_{k+1}$.

Set
\[ \Phi_{n,k}(\lambda) = \int \frac{Q_n(x) d\sigma_k(x)}{\lambda-x}, \qquad k=1,\ldots,p.
\]
In \cite{DL} the authors give a different representation of the polynomials $P_{n,k}$ in terms of minors of the Riemann Hilbert matrix
\[ Y_n(\lambda) := \left(
\begin{array}{cccc}
Q_n(\lambda) & \Phi_{n,1}(\lambda) & \ldots & \Phi_{n,p}(\lambda) \\
Q_{n-1}(\lambda) & \Phi_{n-1,1}(\lambda) & \ldots & \Phi_{n-1,p}(\lambda) \\
\vdots &  \vdots & \ddots & \vdots \\
Q_{n-p}(\lambda) & \Phi_{n-p,1}(\lambda) & \ldots & \Phi_{n-p,p}(\lambda) \\
\end{array}
\right).
\]
For $k=1,\ldots,p-1$, define
\begin{equation} \label{Bnk}
B_{n,k}(\lambda) := \det
\left(
\begin{array}{cccc}
Q_n(\lambda) & \Phi_{n,1}(\lambda) & \ldots & \Phi_{n,k}(\lambda) \\
Q_{n-1}(\lambda) & \Phi_{n-1,1}(\lambda) & \ldots & \Phi_{n-1,k}(\lambda) \\
\vdots &  \vdots & \ddots & \vdots \\
Q_{n-k}(\lambda) & \Phi_{n-k,1}(\lambda) & \ldots & \Phi_{n-k,k}(\lambda) \\
\end{array}
\right).
\end{equation}
Then, see \cite[Proposition 2.6]{DL}
\begin{equation} \label{BP} B_{n,k}(\lambda) = (-1)^{n(k+1)-\frac{k(k+1)}{2}} c_k P_{n,k}(\lambda),\qquad k=1,\ldots,p-1,
\end{equation}
where
\[ c_k = (-1)^k \displaystyle{\left(\int d\mu_1(x)\right)\left(\int Q_1(x)d\mu_2(x)\right)\cdots\left(\int Q_{k-1}(x)d\mu_k(x)\right)}.
\]

The functions $\Psi_{n,k}$ and $\Phi_{n,k}$ are connected. In particular,  $\Psi_{n,1}=\Phi_{n,1}$.  Consequently,
\begin{equation} \label{Bnk2}
B_{n,1}(\lambda) := \det
\left(
\begin{array}{cc}
Q_n(\lambda) & \Psi_{n,1}(\lambda)   \\
Q_{n-1}(\lambda) & \Psi_{n-1,1}(\lambda)   \\
\end{array}
\right).
\end{equation}
Due to \eqref{BP} the zeros of $B_{n,1}$ and $P_{n,1}$ coincide.

\begin{pro} \label{interlacing} $P_{n,1}$ has exactly $N_{n,1} -1$ zeros, they are simple and interlace the zeros of $\Psi_{n,1}$ on $\Gamma_{2}$.
\end{pro}

{\bf Proof.}
Due to \eqref{Bnk2}
\begin{equation}
\label{order2}
 B_{n,1}(z) = Q_n(z) \Psi_{n-1,1}(z) - Q_{n-1}(z) \Psi_{n,1}(z).
\end{equation}
Fix two consecutive zeros of $\Psi_{n,1}$ on $\Gamma_2$, say $x_j, j=1,2$. Then
\[ B_{n,1}(x_j) = Q_n(x_j) \Psi_{n-1,1}(x_j), \qquad j=1,2.
\]
Between $x_1$ an $x_2$ there is exactly one simple zero of $\Psi_{n-1,1}$, because its zeros and those of $\Psi_n$ interlace. The proof of this statement may be carried out in the same fashion as the interlacing property of the zeros of the polynomials $Q_n$ (see, for example, \cite[Theorem 2.1]{AptLopRoc}). On the other hand, $Q_n$ has constant sign on $\Gamma_2$ so
\[ \mbox{sign} (Q_n(x_1) \Psi_{n-1,1}(x_1)) \neq \mbox{sign} (Q_n(x_2) \Psi_{n-1,1}(x_2)).
\]
We conclude that $B_{n,2}$ has an intermediate zero. Taking into consideration the order at infinity of the functions $\Psi_{n,1}$ from \eqref{order2} one sees that $\deg B_{n,1} \leq N_{n,1}-1$; therefore, $\deg B_{n,1} = N_{n,1}-1$, its zeros are simple, and interlace those of $\Psi_{n,1}$ on $\Gamma_2$. \hfill $\Box$

Define
\[\widetilde{s}_{n,k} = \frac{1}{n} \sum_{\Psi_{n,k}(\lambda)=0} \delta_{\lambda}, \qquad k=1,\ldots,p-1,\]
where in the sum each $\lambda$ is counted according to its multiplicity as a zero of $\Psi_{n,k}$ in $\mathbb{C} \setminus \Gamma_k$.
Using the interlacing property proved in Proposition \ref{interlacing}, the fact that $s_2$ is absolutely continuous with respect to the Lebesgue measure on $\Gamma_2$ (thus, it has no mass points), and \eqref{weak},
it is not hard to prove that \eqref{weak} holds for $k=1$ if we replace   $s_{n,1}$ by $\widetilde{s}_{n,1}$. In particular, a property similar to \eqref{zeros2} is verified by the zeros of the functions $\Psi_{n,1}$ on $\Gamma_2$.

\begin{rem} For $k=2,\ldots,p-1$, it would be interesting to prove:
\begin{itemize}
\item  \eqref{weak}  replacing   $s_{n,k}$ by $\widetilde{s}_{n,k}$.
\item a property similar to \eqref{zeros2} for the zeros of the functions $\Psi_{n,k}$ on $\Gamma_{k+1}$.
\item an interlacing property between the zeros of $Q_{n,k}$ and $P_{n,k}$.
\end{itemize}
\end{rem}

\section{An Example}
Here we consider a cubic polynomial, where $q$ has two real negative roots and
a real positive one. A special case of this example was considered by Coussement
et. al \cite{CCV}, and Duits and Kuijlaars \cite{kn:DK}. We begin with
\begin{equation}\label{p}
p(z) := z(a(z)-\lambda)=a_2 z^3+a_1 z^2-\lambda z +1.
\end{equation}

In this case
\begin{equation}\label{qpol}
q(z)=z^3-a_1z-2a_2.
\end{equation}
Let $x_1 < x_2 < 0$. If we set
\begin{equation}\label{coeffp}
a_1=x_1^2+x_2^2+x_1 x_2, \qquad  a_2=-\frac{x_1 x_2}{2}(x_1+x_2),
\end{equation}
then $q$
has as roots $x_1, x_2$ and $-x_1-x_2$. Also
$r$ is given by
\begin{equation}\label{req}
r(z)=z+a_1/z+a_2/z^2.
\end{equation}
The images of the roots of $q$ under $r$ are $\lambda _1=r(x_1)=\frac{4x_1^2+x_2^2+x_1 x_2}{2x_1}$,  $\lambda_2=r(x_2)=\frac{4x_2^2+x_1^2+x_1 x_2}{2x_2}$ and  $\lambda_3=r(-x_1-x_2)=-\frac{4x_1^2+4x_2^2+7x_1 x_2}{2(x_1+x_2)}$.

We look for the root of $p(z)$  that tends to zero as $\lambda$ tends to infinity which is
\begin{equation}\label{root1}
z_0(\lambda)=\frac{1}{6 a_2}\left(u(\lambda)+\frac{l(\lambda)}{u(\lambda)}\right)-\frac{a_1}{3a_2},
\end{equation}
where
\begin{equation}\label{uz}
u(\lambda)=(-k_1(\lambda)+12\sqrt{3}a_2\sqrt{k_2(\lambda)})^{1/3},
\end{equation}
with
\begin{equation}\label{k1z}
k_1(\lambda)=36a_1 a_2 \lambda+ 8 a_1^3+108 a_2^2,
\end{equation}
\begin{equation}\label{kz2}
k_2(\lambda)= -4 a_2 (\lambda-r(x_1))(\lambda-r(x_2))(\lambda-r(-x_1-x_2)),
\end{equation}
and
\begin{equation}\label{lz}
l(\lambda)=4(3a_2 \lambda+a_1^2).
\end{equation}
The principal branch of the logarithm is the one used above. It is not difficult to show that $z_0(\lambda)=1/\lambda +\mathcal{O}(1/\lambda^2).$
Since $x_1<x_2<0$, and $a_2>0$ the cuts are $\Gamma_1=[r(x_1),r(-x_2-x_1)]$ and $\Gamma_2=[- \infty,r(x_2)]$.

The measure
$\rho_1$ will be supported on $\Gamma_1$. On this interval we see
that $k_2(\lambda)\ge0$ and since $k_1(r(x_1))=-(x_2+2x_1)^3(x_2-x_1)^3$ and
$k_1(r(-x_1-x_2))=(x_1+2x_2)^3(2x_1+x_2)^3$, $k_1(\lambda)$ is also positive. Furthermore, from the Cardano formulas we find that
$$
\frac{u(\lambda)}{6a_2}=\left(-\frac{k_1(\lambda)}{(6a_2)^3}+\sqrt{\left(\frac{k_1(\lambda)}{(6a_2)^3}\right)^2-\left(\frac{l(\lambda)}{36a_2^2}\right)^3}\right)^{1/3},
$$
so $u^3(\lambda)$ is negative on this interval. Likewise, $l(\lambda)$ is positive
on $\Gamma_1$.

Examination of the phase shows that for small positive
imaginary $\lambda$ the imaginary part of $u^3$ is positive. Thus, we find that
$$
\frac{u(\lambda)}{6a_2}=(1+\sqrt{3}i)\left(\frac{k_1(\lambda)}{(6a_2)^3}-\sqrt{\left(\frac{k_1(\lambda)}{(6a_2)^3}\right)^2-\left(\frac{l(\lambda)}{36a_2^2}\right)^3}\right)^{1/3}
$$ for $\lambda \in \Gamma_1$,
and
$$
\frac{l(\lambda)}{6a_2 u(\lambda)}=(1-\sqrt{3}i)\left(\frac{k_1(\lambda)}{(6a_2)^3}+\sqrt{\left(\frac{k_1(\lambda)}{(6a_2)^3}\right)^2-\left(\frac{l(\lambda)}{36a_2^2}\right)^3}\right)^{1/3}.
$$
This implies that on $\Gamma_1$,
$$
\rho_1'(x)=-\frac{\sqrt{3}}{\pi}\times
$$
$$\left(\left(\frac{k_1(x)}{(6a_2)^3}+\sqrt{\left(\frac{k_1(x)}{(6a_2)^3}\right)^2-\left(\frac{l(x)}{36a_2^2}\right)^3}\right)^{1/3}-
\left(\frac{k_1(x)}{(6a_2)^3}-\sqrt{\left(\frac{k_1(x)}{(6a_2)^3}\right)^2-\left(\frac{l(x)}{36a_2^2}\right)^3}\right)^{1/3}\right).
$$

In order to compute $\rho_2'$ on $\Gamma_2$, we use that
$$
z_1(\lambda)=\frac{a_1}{3a_2}-\frac{1}{6a_2}\left(u(\lambda)+\frac{l(\lambda)}{u(\lambda)}\right)+i\frac{\sqrt{3}}{12a_2}
\left(u(\lambda)-\frac{l(\lambda)}{u(\lambda)}\right).
$$
The fact that $ u^3(x)$ is positive for $x \in \Gamma_2$ shows that
\[
\rho_1'(x)= \frac{\sqrt{3}}{2\pi(x-r(x_2))} \times \]
\[\left(\left(\frac{-k_1(x)}{(6a_2)^3}+\sqrt{\left(\frac{k_1(x)}{(6a_2)^3}\right)^2-\left(\frac{l(x)}{36a_2^2}\right)^3}\right)^{1/3}+
\left(\frac{k_1(x)}{(6a_2)^3}+\sqrt{\left(\frac{k_1(x)}{(6a_2)^3}\right)^2-\left(\frac{l(x)}{36a_2^2}\right)^3}\right)^{1/3}\right).
\]

\end{document}